\documentclass[dvips]{arxbj}
\usepackage{upgreek,mathbh}

\aid{0}
\volume{16}
\issue{3}
\pubyear{2010}
\firstpage{825}
\lastpage{857}
\doi{10.3150/09-BEJ234}

\makeatletter
\def\epsilon{\varepsilon}

\newproclaim{kk1}{KK1}
\newproclaim{kk2}{KK2}
\newproclaim{kk3}{KK3}
\newremark{remarks}{Remarks}
\newremark{remark}{Remark}
\newproclaim{ks1}{KS1}
\newproclaim{ks2}{KS2}
\newtheorem{theorem}{Theorem}
\newtheorem{proposition}{Proposition}
\newtheorem{corollary}{Corollary}
\newtheorem{lemma}{Lemma}

\makeatother

\begin{document}
\begin{frontmatter}

\title{A self-similar process arising from a random walk with random
environment in~random~scenery}
\runtitle{A self-similar process}

\begin{aug}
\author[a]{\fnms{Brice} \snm{Franke}\corref{}\thanksref{a}\ead[label=e1]{Brice.Franke@rub.de}}
\and
\author[b]{\fnms{Tatsuhiko} \snm{Saigo}\thanksref{b}\ead[label=e2]{saigo@math.keio.ac.jp}}
\runauthor{B. Franke and T. Saigo}
\address[a]{Fakult\"at f\"ur Mathematik,
Ruhr-Universit\"at Bochum, Universit\"atsstr. 150, 44780 Bochum,
Germany. \mbox{\printead{e1}}}
\address[b]{Department of Mathematics, Keio University 3-14-1 Hiyoshi,
Kouhoku-ku,
Yokohama-shi City,\\ Kanagawa-ken Prefecture, 223-8522, Japan. \printead{e2}}
\end{aug}

\received{\smonth{1} \syear{2009}}
\revised{\smonth{5} \syear{2009}}

%
\begin{abstract}
In this article, we merge celebrated results of Kesten and Spitzer
[\textit{Z. Wahrsch. Verw. Gebiete} \textbf{50} (1979) 5--25]
and Kawazu and Kesten [\textit{J. Stat. Phys.} \textbf{37} (1984)
561--575]. A random walk performs a motion in an i.i.d. environment and
observes an i.i.d. scenery along its path. We assume that the scenery
is in the domain of attraction
of a stable distribution and prove that the resulting observations
satisfy a limit theorem.
The resulting limit process is a self-similar stochastic process with
non-trivial
dependencies.
\end{abstract}

%

\end{frontmatter}

\section{Introduction}\label{sec1}

The following model for a random walk in a random environment can be
found in the physics literature; see Anshelevic and Vologodskii (\citeyear{AnsVol1981}),
Alexander \textit{et al.} (\citeyear{Aleetal1981}), Kawazu and Kesten (\citeyear{KawKes1984}).
Let $ \{\lambda_j;j\in\mathbb{Z}\} $ be a family of positive
i.i.d. random variables and $ \mathcal{A} $ the $ \sigma$-algebra
generated by those random variables.
Let $ \{X(t);t\geq0\} $ be a continuous-time random walk on $ \mathbb
{Z} $ having the following asymptotic
transition rates for $ h\rightarrow0$:
\begin{eqnarray} \label{Formel1}
\mathbb{P} \bigl(X(t+h)=j+1|X(t)=j,\mathcal{A}\bigr) &=& \lambda_jh+\mathrm{o}(h),\\
\mathbb{P} \bigl(X(t+h)=j-1|X(t)=j,\mathcal{A}\bigr) &=& \lambda_{j-1}h+\mathrm{o}(h),\\
\mathbb{P} \bigl(X(t+h)=j|X(t)=j,\mathcal{A}\bigr) &=& 1-(\lambda_j+\lambda_{j-1})h+\mathrm{o}(h).
\end{eqnarray}
In other words, the process $ \{X(t);t\geq0\} $ is a birth--death
process with possibly negative
population size, where, for a population with $ j $ individuals, birth
occurs at rate $ \lambda_j $ and
death at rate $ \lambda_{j-1} $.
We will assume that the process $ \{X(t);t\geq0\} $ starts at zero at
time zero.
The resulting process is symmetric, in the sense that the permeability
of the edge connecting
the vertices $ j $ and $ j+1 $ does not depend on the direction of the
motion. This physical background
motivates the name `random environment' for the sequence $ \{\lambda
_j;j\in\mathbb{Z}\} $.\ In what follows, we
denote the distribution of the random environment on the sequence space
by $ P_\lambda$.
The following convergence results are described in Kawazu and Kesten
(\citeyear{KawKes1984}).

\begin{kk1*}
If $ c:=\mathbb{E} [\lambda_0^{-1}]<\infty$, then
for $
P_\lambda$-almost all environments,
the distributions (after conditioning on the environment) of the
processes
\[
X_n(t):=\frac{1}{n}X(n^2t),\qquad t\geq0,
\]
converge weakly with respect to the Skorohod topology toward the
distribution
of the process
$ \{c^{-1/2}B(t);t\geq0\} $, where $ \{B(t);t\geq0\} $ is standard
Brownian motion on $ \mathbb{R} $.
\end{kk1*}

(See also Papanicolaou and Varadhan (\citeyear{PapVar1981}) for some related results.)

\begin{kk2*}
If there exists a slowly varying function $ L_1 $ such that
\[
\frac{1}{nL_1(n)}\sum_{j=1}^n\frac{1}{\lambda_j}\longrightarrow1\qquad
\mbox{in probability},
\]
then the distributions of the processes
\[
X_n(t):=\frac{1}{n}X(n^2L_1(n)t)
\]
converge weakly with respect to the Skorohod topology toward the
distribution of standard Brownian
motion.
\end{kk2*}

\begin{kk3*}
If there exists a slowly varying function $ L_2 $ such
that the sequence of random variables
\[
R_n:=\frac{1}{n^{1/\alpha}L_2(n)}\sum_{j=1}^n\frac{1}{\lambda_j}
\]
converges in distribution toward a one-sided stable distribution $
\vartheta_\alpha$ with
index $ \alpha\in(0,1) $, then the distributions of the processes
\[
X_n(t):=\frac{1}{n}X\bigl(n^{(1+\alpha)/\alpha}L_2(n)t\bigr)
\]
converge weakly with respect to the Skorohod topology toward the
distribution of a continuous
self-similar process $ \{X_\ast(t);t\geq0\} $ with scaling exponent $
\eta=\frac{\alpha}{\alpha+1} $.
\end{kk3*}

\begin{remarks*}
(1) In the next section, we will give a representation for the process $
X_\ast$ in terms of a standard Brownian motion and a
stable subordinator associated with the measure $ \vartheta_\alpha$.

(2) We note that the results from Kawazu and Kesten (\citeyear{KawKes1984}) are
generalized in Kawazu (\citeyear{Kaw1989}).
\end{remarks*}

He considered random walks in random environments defined by the
following transition asymptotics:
\begin{eqnarray*}
\mathbb{P} \bigl(X(t+h)=j+1|X(t)=j,\mathcal{A}\bigr) &=& (\lambda_j/\eta
_j)h+\mathrm{o}(h),\\
\mathbb{P} \bigl(X(t+h)=j-1|X(t)=j,\mathcal{A}\bigr) &=& (\lambda_{j-1}/\eta
_j)h+\mathrm{o}(h),\\
\mathbb{P} \bigl(X(t+h)=j|X(t)=j,\mathcal{A}\bigr) &=& 1-\bigl((\lambda_j+\lambda
_{j-1})/\eta_j\bigr)h+\mathrm{o}(h),
\end{eqnarray*}
where $ \{\eta_j,j\in\mathbb{N}\} $ is an i.i.d. family of positive
random variables satisfying suitable
assumptions.
Similarly to the situation studied in Kawazu and Kesten (\citeyear{KawKes1984}), the
resulting random walks converge toward
appropriate continuous processes after scaling.

In Kesten and Spitzer (\citeyear{KesSpi1979}), new classes of continuous self-similar
processes are described.
Moreover, it was proven therein that those processes are weak limits of
random walks in random scenery.
Those random walks are defined as follows.

Let $ \{\xi(x);x\in\mathbb{Z}\} $ and $ \{Z_i;i\in\mathbb{N}\} $
be two independent families of i.i.d. random variables, where the
random variables $ Z_i $ are assumed to be $ \mathbb{Z} $-valued.
One can think of the sequence $ \{Z_i;i\in\mathbb{N}\} $ as
increments of a classical
$ \mathbb{Z} $-valued random walk $ S_k:=\sum_{i=1}^kZ_i $.
The stationary sequence $ \{\xi(S_k);k\in\mathbb{N}\} $ has some
non-trivial long-range dependencies if
the underlying random walk $ \{S_k;k\in\mathbb{N}\} $ is recurrent.
This is the case, for example, if $ Z_1 $ is in the domain of
attraction of an $ \alpha$-stable distribution
with $ \alpha\in(1,2] $. The random sequence $ D(n):=\sum_{k=1}^n\xi
(S_k) $ is called a
\emph{random walk in random scenery}.
In Kesten and Spitzer (\citeyear{KesSpi1979}), the following convergence result
was proven for those processes.

\begin{ks1*}
If $ \xi(0) $ is in the domain of attraction of a $
\beta$-stable distribution
with $ \beta\in(0,2] $ and if $ Z_1 $ is in the domain of attraction
of an $ \alpha$-stable distribution with
$ \alpha\in(0,1) $, then the distributions of the processes
\[
D_n(t):=n^{-1/\beta}\sum_{k=1}^{\lfloor nt\rfloor}\xi(S_k)
\]
\it converge weakly with respect to the Skorohod topology toward $
\beta$-stable L\'{e}vy motion.
\end{ks1*}

(See also Spitzer (\citeyear{Spi1976}) for a special case.)

\begin{ks2*}
If $ \xi(0) $ is in the domain of attraction of a $
\beta$-stable distribution
with $ \beta\in(0,2] $ and if $ Z_1 $ is in the domain of attraction
of an $ \alpha$-stable distribution with
$ \alpha\in(1,2] $, then the distributions of the processes
\[
D_n(t):=n^{-\delta}\sum_{k=1}^{\lfloor nt\rfloor}\xi(S_k)
\]
converge weakly with respect to the Skorohod topology toward a
continuous self-similar process~$ D_\ast$ with scaling exponent $ \delta=1-\frac{1}{\alpha}+\frac
{1}{\alpha\beta}$.
\end{ks2*}

\begin{remark*}
The statement in KS1 corresponds to the transient case
and is
not difficult to prove since, in that case, the sequence
$ \{\xi(S_k);k\in\mathbb{N}\} $ has only weak dependencies.
This is the reason why one obtains $ \beta$-stable L\'evy noise in the
limit. We also mention that the case $ \beta=1 $ is still open.
\end{remark*}

\begin{remark*}
There exist various generalizations of the results of Kesten and Spitzer (\citeyear{KesSpi1979}).
We will only mention Shieh (\citeyear{Shi1995}), where the limiting process is
generalized to higher dimensions,
Lang and Nguyen (\citeyear{LanNgu1983}), which deals with multidimensional random walks
and some special
random scenery, Maejima (\citeyear{Mae1996}), where the random scenery belongs to the
domain of attraction of
an operator-stable distribution, Arai (\citeyear{Ara2001}), where the random scenery
belongs to the domain of partial
attraction of a semi-stable distribution, and Saigo and Takahashi
(\citeyear{SaiTak2005}), where the random scenery
and the random walk belong to the domain of partial attraction of
semi-stable and operator semi-stable distributions.
\end{remark*}

In this article, we investigate whether it is possible to substitute
the classical
random walk in the result of Kesten and Spitzer (\citeyear{KesSpi1979}) by the random
walk in random environment
which was introduced in Kawazu and Kesten (\citeyear{KawKes1984}).
We will restrict our attention to the result KK3 since this is the case
where a new type
of self-similar process arises at the end. For simplicity and in order
to avoid complicating notation,
we will assume that the slowly varying function $ L_2 $ which appears
in KK3 is constant and equal
to one. The general case involving non-constant $ L_2 $ can be treated
in a similar way.

We now fix a probability space $ (\Omega,\mathcal{F},\mathbb{P} ) $
which is
sufficiently large to support a
family of i.i.d. random variables $ \{\lambda_j;j\in\mathbb{Z}\} $,
a birth--death process
$ \{X(t);t\geq0\} $ with asymptotic transition rates given by
equations (1)--(3) and a family of
i.i.d. random variables $ \{\xi(k),k\in\mathbb{Z}\} $.

We assume that the families $ \{\xi(k),k\in\mathbb{Z}\} $ and $ \{
X(t);t\geq0\} $ are independent and
that $ t\mapsto X(t) $ is cadlag $ \mathbb{P} $-almost surely.

Further, we assume
that $ \lambda^{-1}_1 $ is in the domain of normal attraction of a
one-sided $ \alpha$-stable distribution
$ \vartheta_\alpha$ with $ \alpha\in(0,1) $.

Moreover, we assume that $ \xi(0) $ is in the domain of normal
attraction of a $ \beta$-stable distribution
$ \vartheta_\beta$ with $ \beta\in(0,2] $. Its characteristic
function is given by
\[
\psi(\theta)=\exp\bigl(-|\theta|^\beta\bigl(A_1+\mathrm{i}A_2\operatorname{sgn}(\theta)\bigr)\bigr) ,
\]
where $ 0<A_1<\infty$ and $ |A_1^{-1}A_2|\leq\tan(\uppi\beta/2) $.
For $ \beta>1 $, it follows from those assumptions that $ \mathbb{E}
[\xi(0)]=0 $.

For $ \beta=1 $, we make the further assumption that there exists a $
K>0 $ such that
\[
\bigl|\mathbb{E} \bigl[\xi(0)\mathbh{1}_{[-\rho,\rho]}(\xi(0)) \bigr]
\bigr|\leq K\qquad
\mbox{for all } \rho>0 .
\]
We can now define the following continuous-time version of the random
walk in random scenery:
\[
\Xi(t):=\int_0^t\xi(X(s))\,\mathrm{d}s .
\]
In the following, we will use the space
\[
D[0,\infty):= \{\gamma\dvtx[0,\infty)\rightarrow\mathbb{R}\dvtx\gamma
\mbox{ is cadlag} \}
\]
with the Skorohod topology.
We will prove the following theorem.

\begin{theorem} \label{MT}
For $ \kappa:=\frac{1}{\alpha}+\frac{1}{\beta} $ and $
k_n:=n^{(1+\alpha)/\alpha}$,
the distributions of the processes
\[
\Xi_n(t):=n^{-\kappa}\int_0^{k_nt}\xi(X(s))\,\mathrm{d}s
\]
converge weakly with respect to the Skorohod topology toward the
distribution of a
self-similar stochastic process $ \{\Xi_\ast(t);t\geq0\} $ with
scaling exponent
$ \mu=1-\frac{\alpha}{\alpha+1}+\frac{\alpha}{(\alpha+1)\beta} $.
\end{theorem}

\begin{remark*}
The stochastic process $ \{\Xi_\ast(t);t\geq0\} $ can
be constructed as follows. Let $ Z_+ $ and $ Z_- $ be two independent
copies of the $ \beta$-stable
L\'{e}vy process which can be associated with the characteristic function
\[
\psi(\theta)=\exp\bigl(-|\theta|^\beta\bigl(A_1+\mathrm{i}A_2\operatorname{sgn}(\theta)\bigr) \bigr) .
\]
Further, let $ \{L_\ast(\tau,x);\tau\geq0,x\in\mathbb{R}\} $ be
the local time of the stochastic process
$ \{X_\ast(\tau);\tau\geq0\} $; that is, the random variable $
L_\ast(\tau,x) $ is the derivative with respect to
$ x $ of the occupation time
\[
\Gamma_\ast(\tau,(-\infty,x]):=\int_0^\tau\mathbh{1}_{(-\infty,x]}(X_\ast(\sigma))\,\mathrm{d}\sigma.
\]
We will see in the next section that the local time exists for all but
a countable number of points
$ x\in\mathbb{R} $.
Moreover, for all $ \tau\geq0 $, the processes
\[
\{L_\ast(\tau,x-);x\geq0\} \quad \mbox{and} \quad  \{L_\ast(\tau
,-(x-));x\geq0\}
\]
are predictable with respect to the natural filtrations of $ Z_+ $
(resp., $ Z_-$).
The following integral representation of the process $ \Xi_\ast$ can
be given:
\[
\Xi_\ast(\tau):=\int_0^\infty L_\ast(\tau,x-)\,\mathrm{d}Z_+(x)+\int
_0^\infty L_\ast(\tau,-(x-))\,\mathrm{d}Z_-(x) .
\]
\end{remark*}

\section{The convergence of the birth--death process}

The goal of this section is to prove Corollary \ref{PrinceKor}, which
is the main ingredient needed to show that
the finite-dimensional distributions of $ \Xi_n $ converge toward the
finite-dimensional distributions of
$ \Xi_\ast$.
This corollary contains a statement on the weak convergence of certain
functionals of the occupation times
of the rescaled processes
$ X_n $. A result corresponding to Corollary \ref{PrinceKor} is also
proved in Kesten and Spitzer (\citeyear{KesSpi1979});
however, we have to adopt a totally different approach since we do not
have such precise information on the
potential theory related to the random walk~$ X $.
Instead, we will understand the occupation times of $ X_n $ and prove
that they converge in an appropriate
sense toward the local time of the limit process $ X_\ast$.

We describe some of the main arguments from the proof in Kawazu and Kesten (\citeyear{KawKes1984}) for the convergence
of the processes
\[
X_n(t):=\frac{1}{n}X\bigl(n^{(1+\alpha)/{\alpha}}t\bigr)
\]
toward the self-similar process $ X_\ast$ defined in Kawazu and Kesten (\citeyear{KawKes1984}).
We can enlarge our underlying probability space $ (\Omega,\mathcal{F},\mathbb{P}
) $
in such a way that it contains a standard Brownian motion $ \{
B(t);t\geq0\} $ and a cadlag version of the
stable L\'{e}vy subordinator $ \{W(x);x\in\mathbb{R}\} $ which can be
associated with the one-sided
$ \alpha$-stable distribution $ \vartheta_\alpha$.

Furthermore, we assume that $ \{B(t);t\geq0\} $, $ \{W(x);x\in\mathbb
{R}\} $, $ \{X(t);t\geq0\} $ and
$ \{\xi(n);n\in\mathbb{Z}\} $ are independent.
Moreover, we assume that $ W(0)=0 $ and $ B(0)=0 $ hold $ \mathbb{P} $-almost
surely.

In the future, we will denote by $ \{L(t,x);t\geq0,x\in\mathbb{R}\}
$ the local time of the Brownian motion
$ \{B(t);t\geq0\} $.
The process
\[
V_\ast(t):=\int_\mathbb{R}L(t,W(x))\,\mathrm{d}x
\]
is non-decreasing $ \mathbb{P} $-almost surely.
Therefore, we can define the following pseudo-inverse:
\[
W^{-1}(y):=\inf\{x\in\mathbb{R};W(x)>y\} \quad \mbox{and} \quad
V_\ast^{-1}(\tau):=\inf\{t\geq0;V_\ast(t)>\tau\} .
\]
In Kawazu and Kesten (\citeyear{KawKes1984}), the following representation for the
self-similar process $ X_\ast$ is given:
\[
X_\ast(\tau):=W^{-1}(B(V_\ast^{-1}(\tau))).
\]
We now sketch the main arguments from the proof in Kawazu and Kesten (\citeyear{KawKes1984}).
We will need some of those ideas in our proof of the convergence of $
\Xi_n $ toward $ \Xi_\ast$.
Their approach is based on the natural scale of the birth--death
process. One defines
\[
S(j):= \cases{\displaystyle
\sum_{k=0}^{j-1}\lambda_k^{-1} &\quad for $ j>0$,\vspace*{2pt}\cr
0 & \quad for $j=0$,\cr
\displaystyle-\sum_{k=j}^{-1}\lambda_k^{-1} &\quad for $ j<0$.
}
\]
This implies that conditioned on $ \mathcal{A}:=\{\lambda_j;j\in\mathbb
{Z}\}, $ the process $ S(X(t)) $ is on
natural scale (see Kawazu and Kesten (\citeyear{KawKes1984}), page 565).
This means that for all $ a,b,x\in\mathbb{R} $ with $ a<x<b $, one has
\[
\mathbb{P} \bigl(S(X(t)) \mbox{ hits } \{a,b\} \mbox{ first at } a
\mid S(X(0))=x,\mathcal{A}\bigr)=\frac{b-x}{b-a}.
\]
It is then possible to represent the process $ S(X(t)) $ as the time
change of standard Brownian motion
$ \{B(t);t\geq0\} $ as follows.

One defines $ m(\mathrm{d}x):=\sum_{i\in\mathbb{Z}}\delta_{S(i)}(\mathrm{d}x) $ and
\[
V(t):=\int_\mathbb{R} L(t,x)m(\mathrm{d}x)=\sum_{i\in\mathbb{Z}}L(t,S(i)) ,
\]
where $ \{L(t,x);t\geq0,x\in\mathbb{R}\} $ is again the local time
of the standard Brownian motion $ B $.
One can see that $ \{B(V^{-1}(t));t\geq0\} $ and $ \{S(X(t));t\geq0\}
$ are both cadlag and have the same
distribution (see Kawazu and Kesten (\citeyear{KawKes1984}), page 566).

One then has to scale the above constructions.
\[
S_n(x):=n^{-1/\alpha}S(\lfloor nx\rfloor),\qquad n\in\mathbb{N},
x\in\mathbb{R} ,
\]
where, for a positive real number $ x $, we denote by $ \lfloor
x\rfloor$ its integer part.
It follows from the assumptions on the environment $ \{\lambda_j;j\in
\mathbb{Z}\} $ that for
$ n\rightarrow\infty$, the processes $ \{S_n(x);x\in\mathbb{R}\} $
converge in distribution toward an
$ \alpha$-stable L\'{e}vy process $ \{W(x);x\in\mathbb{R}\} $.
Moreover, the process $ W $ is strictly increasing $ \mathbb{P}
$-almost surely since
$ \vartheta_\alpha$ is a one-sided stable distribution and $ \alpha
\in(0,1) $.
By a method given in Skorohod (\citeyear{Sko1956}) and Dudley (\citeyear{Dud1968}), it is possible
to construct a suitable probability
space $ (\tilde{\Omega},\tilde\mathcal{F},\tilde{\mathbb{P} }) $ with suitable
$ D $-valued random variables
$ \tilde{S}_n $ and $ \tilde{W} $ having the properties that $ \tilde
{S}_n $ converges toward $ \tilde{W} $
almost surely with respect to $ \tilde{\mathbb{P} } $ and that $
\tilde{S}_n
$ and $ \tilde{W} $ have the same
distributions as $ S_n $ (resp., $ W $) (see Kawazu and Kesten (\citeyear{KawKes1984}),
page 567).
One then defines
\[
\tilde{V}_n(t):=\int_\mathbb{R} L(t,x)\tilde{m}_n(\mathrm{d}x) \quad
\mbox{and}\quad
\tilde{V}_\ast(t):=\int_\mathbb{R} L(t,x)\tilde{m}_\ast(\mathrm{d}x)
\]
with
\[
\int_\mathbb{R} f(x)\tilde{m}_n(\mathrm{d}x):=\int_{\mathbb{R}}f(\tilde
{S}_n(x))\,\mathrm{d}x\quad \mbox{and}\quad
\int_\mathbb{R} f(x)\tilde{m}_\ast(\mathrm{d}x):=\int_{\mathbb{R}}f(\tilde
{W}(x))\,\mathrm{d}x
\]
for all measurable $ f\geq0 $. We then define $ \tilde{S}_n^{-1} $, $
\tilde{W}^{-1} $,
$ \tilde{V}_n^{-1} $ and $ \tilde{V}_\ast^{-1} $ in the same way as
$ W^{-1} $ (resp., $ V_\ast^{-1} $)
above.

In Kawazu and Kesten (\citeyear{KawKes1984}) (see page 568) they prove that
$ \{B(\tilde{V}^{-1}_n(t));t\geq0\} $ converges $ \tilde{\mathbb{P} }
$-almost surely toward
$ \{B(\tilde{V}^{-1}_\ast(t));t\geq0\} $ in the $ J_1 $-topology.
For convenience, we define
\[
\tilde{X}_n(t):=\tilde{S}_n^{-1}(B(\tilde{V}^{-1}_n(t))), \qquad
\tilde{X}_\ast(t):=\tilde{W}^{-1}(B(\tilde{V}^{-1}_\ast(t))) .
\]
We note that the process
$ \{\tilde{X}_n(t);t\geq0\} $ is defined on
$ (\Omega\times\tilde{\Omega},\mathcal{F}\times\tilde\mathcal{F},\mathbb{P}
\times\tilde{\mathbb{P} }) $.
It is proved in Kawazu and Kesten (\citeyear{KawKes1984}) that $ \{\tilde{X}_n(t);t\geq
0\} $ converges toward
$ \{\tilde{X}_\ast(t);t\geq0\} $ with respect
to the $ J_1 $-topology almost surely with respect to $ \mathbb{P}
\times
\tilde{\mathbb{P} } $ (see page 569).

Moreover, for $ B_n(t):=n^{-1/2}B(nt) $ one has that (see Kawazu and Kesten (\citeyear{KawKes1984}), page 572)
\[
|X_n(t)-S_n^{-1}(B_n(V_n^{-1}(t)))|\leq1/n
\]
and
\[
\{S_n^{-1}(B_n(V^{-1}_n(t)));t\geq0 \}\stackrel{\mathcal{D}}{=}
\{\tilde{S}^{-1}_n(B(\tilde{V}_n^{-1}(t)));t\geq0 \}
= \{\tilde{X}_n(t);t\geq0 \} .
\]
If we define $ \hat{X}_n(t):=S_n^{-1}(B_n(V^{-1}_n(t))) $, then the
previous observations imply that
both processes $ \{X_n(t);t\geq0\} $ and $ \{\hat{X}_n(t);t\geq0\} $
converge in distribution toward
$ \{\tilde{X}_\ast(t);t\geq0\} $, which has the same distribution as
$ \{X_\ast(t);t\geq0\} $.

In the rest of this section, we analyze the distributional behavior of
the occupation times for the process
$ X_n $ (see Proposition \ref{PrinceProp}). In order to obtain this
result, we prove an analogous result for
the process $ \tilde{X}_n $ (see Lemma \ref{PrinceLem}), which can be
reduced to
Proposition \ref{CardTowardMeasureProp}. The advantage of this detour
is that we can prove almost sure
convergence for the occupation times of the process $ \tilde{X}_n $
toward the local time of
$ \tilde{X}_\ast$ (see Proposition \ref{OkTimeLokTimeConvProp}).
This result is based on the fact that
we have explicit formulas for the occupation times of $ \tilde{X}_n $
and the local time of
$ \tilde{X}_\ast$ (see Proposition \ref{OkTimePro} and Corollary
\ref{LokTimeKor1}).
The explicit expression of the occupation time of $ \tilde{X}_n $ and
the local time of $ \tilde{X}_\ast$
reveals that in order to prove Proposition \ref
{OkTimeLokTimeConvProp}, it is sufficient to prove the
almost sure convergence of $ \tilde{S}_n $ and $ \tilde{V}_n^{-1} $
toward $ \tilde{W}_\ast$ (resp., $ \tilde{V}^{-1}_\ast$). The
convergence of $ \tilde{S}_n $ toward $ \tilde{W}_\ast$ holds by
construction. The convergence of $ \tilde{V}_n $ toward $ \tilde
{V}_\ast$ is obtained in Lemma
\ref{TimeChangeConvLem} and then used to obtain the convergence of $
\tilde{V}_n^{-1} $ toward
$ \tilde{V}^{-1}_\ast$ in Lemma \ref{InversTimeChangeConvLem}.

\subsection{The local times of $ X_\ast$ and $ \tilde{X}_\ast$}

We define the time that the processes $ \tilde{X}_\ast$ and $ X_\ast
$ spend in the measurable set $ A $
until time $ \tau$ as
\[
\Gamma_\ast(\tau,A):=\int_0^\tau\mathbh{1}_{A}(X_\ast
(\sigma))\,\mathrm{d}\sigma\qquad \biggl(\mbox{resp.},\
\tilde{\Gamma}_\ast(\tau,A):=\int_0^\tau\mathbh{1}_{A}(\tilde{X}_\ast(\sigma))\,\mathrm{d}\sigma\biggr).
\]
We denote by $ \{L_\ast(\tau,x);\tau\geq0,x\in\mathbb{R}\} $ and
$ \{\tilde{L}_\ast(\tau,x);\tau\geq0,x\in\mathbb{R}\} $ the
local times of $ X_\ast$ (resp.,
$ \tilde{X}_\ast$) if they exist.
In this subsection, we prove that both local times exist almost surely
and relate them to the
local time $ \{L(t,x);t\geq0,x\in\mathbb{R}\} $ of the underlying
Brownian motion $ \{B(t);t\geq0\} $.

\begin{proposition} \label{LokTimePro}
One has $ \mathbb{P} $-almost surely that for $ \tau\geq0 $ and all
$ x\in
\mathbb{R} $,
\[
\Gamma_\ast(\tau,(-\infty,x))=\int_{-\infty}^xL(V_\ast^{-1}(\tau
),W(y))\,\mathrm{d}y .
\]
Further, $ \mathbb{P} \times\tilde{\mathbb{P} } $-almost surely for
all $ \tau\geq0
$ and all
$ x\in\mathbb{R} $,
\[
\tilde{\Gamma}_\ast(\tau,(-\infty,x))=\int_{-\infty}^xL(\tilde
{V}_\ast^{-1}(\tau),\tilde{W}(y))\,\mathrm{d}y .
\]
\end{proposition}

\begin{pf}
We have $ \mathbb{P} $-almost surely that $ x\mapsto W(x) $ is increasing.
It follows that the set $ \mathcal{N}_1 $ of $ x\in\mathbb{R} $ where $
W $ is not continuous
is countable.
We define the set
\[
\mathcal{N}_2:= \bigl\{x\in\mathbb{R}\dvtx\ell\bigl(\sigma;B(V_\ast^{-1}(\sigma
))=W(x)\bigr)>0 \bigr\} ,
\]
where $ \ell$ denotes the Lebesgue measure on $ \mathbb{R} $.
The set $ \mathcal{N}_2 $ is countable since for $ x_1\neq x_2 $, one has
that the sets
$ \{\sigma;B(V_\ast^{-1}(\sigma))=W(x_1)\} $ and $\{\sigma
;B(V_\ast^{-1}(\sigma))=W(x_2)\} $ are disjoint.
The statement then follows since there cannot be an uncountable number
of disjoint subsets of
$ \mathbb{R} $ with positive Lebesgue measure. Thus the set $ \mathcal{N}:=\mathcal{N}_1\cup\mathcal{N}_2 $
is countable.
Since the function $ x\mapsto\Gamma_\ast(\tau,(-\infty,x)) $ is
increasing and since
\[
x\mapsto\int_{-\infty}^xL(V_\ast^{-1}(\tau),W(y))\,\mathrm{d}y
\]
is continuous, it is sufficient to prove the statement of the
proposition for $ x\in\mathcal{N}^c $.

The fact that $ W $ is increasing and continuous in $ x $ implies the
equivalence of the statement
$ W(x)>y $ with the statement $ \exists z_0<x\dvtx W(z_0)>y $.

The latter statement is then equivalent to the statement
$ W^{-1}(y):=\inf\{z\dvtx W(z)>y\}<x $.

This then implies that
$ \mathbh{1}_{(-\infty,x)}(X_\ast(\sigma))
=\mathbh{1}_{(-\infty,W(x))}(B(V_\ast^{-1}(\sigma))) $.

We also note that $ t\mapsto V(t) $ is continuous and non-decreasing.
This implies that
$ V_\ast\circ V_\ast^{-1}=\operatorname{id}_\mathbb{R} $.

In the following, we want to compute the derivative of the
non-decreasing function
\[
M\dvtx \sigma\mapsto\int_{-\infty}^xL(V_\ast^{-1}(\sigma),W(y))\,\mathrm{d}y .
\]
Since $ W $ is increasing and continuous in $ x $, we have that
$ B(V_\ast^{-1}(\sigma_0))<W(x) $ implies that
\[
\sigma\mapsto\int_x^\infty L(V_\ast^{-1}(\sigma),W(y))\,\mathrm{d}y
\]
is locally constant, say equal to $c_0$,
in a  neighborhood of $ \sigma_0 $.

Thus
\[
\sigma\mapsto\int_{-\infty}^xL(V_\ast^{-1}(\sigma),W(y))\,\mathrm{d}y
= V_\ast(V_\ast^{-1}(\sigma))-c_0=\sigma-c_0
\]
in a neighborhood of $ \sigma_0 $.

Moreover, since $ W $ is increasing and continuous in $ x $, we have that
$ B(V_\ast^{-1}(\sigma_0))>W(x)$ implies
\[
\sigma\mapsto\int_{-\infty}^xL(V_\ast^{-1}(\sigma),W(y))\,\mathrm{d}y
\]
is locally constant  in a  neighborhood of $ \sigma_0 $.

It therefore turns out that
\[
M'(\sigma)= \cases{
1, &\quad if $ B(V_\ast^{-1}(\sigma))<W(x)$,\cr
0, &\quad if $ B(V_\ast^{-1}(\sigma))>W(x)$.
}
\]
Moreover, for all $ \sigma_1,\sigma_2\in\mathbb{R}^+ $ with $
\sigma_1\leq\sigma_2 $, we have that
\[
\int_{-\infty}^xL(V_\ast^{-1}(\sigma_1),W(y))\,\mathrm{d}y
\leq\int_{-\infty}^xL(V_\ast^{-1}(\sigma_2),W(y))\,\mathrm{d}y
\]
and
\[
\int_x^\infty L(V_\ast^{-1}(\sigma_1),W(y))\,\mathrm{d}y
\leq\int_x^\infty L(V_\ast^{-1}(\sigma_2),W(y))\,\mathrm{d}y .
\]
This implies that
\begin{eqnarray*}
&&\int_{-\infty}^xL(V_\ast^{-1}(\sigma_2),W(y))\,\mathrm{d}y
-\int_{-\infty}^xL(V_\ast^{-1}(\sigma_1),W(y))\,\mathrm{d}y\\
&&\quad
\leq V_\ast(V_\ast^{-1}(\sigma_2))-V_\ast(V_\ast^{-1}(\sigma
_1))=\sigma_2-\sigma_1.
\end{eqnarray*}
It follows that
\[
\sigma\mapsto\int_{-\infty}^xL(V_\ast^{-1}(\sigma),W(y))\,\mathrm{d}y
\]
is Lipschitz continuous with Lipschitz constant smaller than one.

Since the set $ \{\sigma\dvtx B(V_\ast^{-1}(\sigma))=W(x)\} $ is a zero
set with respect
to the Lebesgue measure $ \ell$ for all $ x\in\mathcal{N}^c $, it
follows that
\[
\int_0^\tau\mathbh{1}_{(-\infty,x)}(X_\ast(\sigma
))\,\mathrm{d}\sigma=
\int_0^\tau\mathbh{1}_{(-\infty,W(x))}(B(V_\ast
^{-1}(\sigma)))\,\mathrm{d}\sigma=
\int_0^\tau M'(\sigma)\,\mathrm{d}\sigma=M(\tau) .
\]
The second statement is proved in the same way.
\end{pf}

\begin{corollary} \label{LokTimeKor1}
One has $ \mathbb{P} $-almost surely that the local time $ L_\ast
(\tau,x) $ is
defined for all $ \tau\geq0 $ and all $ x $, where $ x\mapsto W(x) $
is continuous.
Further, one has $ \mathbb{P} \times\tilde{\mathbb{P} } $-almost
surely that the
local time $ \tilde{L}_\ast(\tau,x) $
is defined for all $ \tau\geq0 $ and all $ x $, where $ x\mapsto
\tilde{W}(x) $ is continuous.
In those points, one has
\[
L_\ast(\tau,x)=L(V_\ast^{-1}(\tau),W(x)) \qquad  \bigl(\mbox{resp.}, \
\tilde{L}_\ast(\tau,x)=L(\tilde{V}_\ast^{-1}(\tau),\tilde
{W}(x))\bigr) .
\]
\end{corollary}

\begin{pf}
Differentiation in Proposition \ref{LokTimePro} proves
this corollary.
\end{pf}

\subsection{The occupation time of $ \tilde{X}_n $}

For a measurable set $ A\subset\mathbb{R} $, we define
\[
\hat{\Gamma}_n(t,A):=\int_0^t \mathbh{1}_{A}(\hat
{X}_n(\sigma))\,\mathrm{d}\sigma, \qquad
\tilde{\Gamma}_n(t,A):=\int_0^t \mathbh{1}_{A}(\tilde
{X}_n(\sigma))\,\mathrm{d}\sigma
\]
and
\[
\Gamma_n(t,A):=\int_0^t \mathbh{1}_{A}(X_n(\sigma
))\,\mathrm{d}\sigma.
\]
These are the respective times that the processes $ \hat{X}_n $, $
\tilde{X}_n $ and $ X_n $ spend in the set $ A $
until time $ t $.
In this section, we give an explicit expression for the occupation time
of $ \tilde{X}_n $ in terms
of the local time $ \{L(t,x);t\geq0,x\in\mathbb{R}\} $ of the
underlying Brownian motion $ \{B(t);t\geq0\} $.

\begin{proposition}\label{OkTimePro}
One has $ \mathbb{P} \times\tilde{\mathbb{P} } $-almost surely for
all $ \tau\geq0
$ and all $ x\in\mathbb{R} $ that
\[
\tilde{\Gamma}_n(\tau,\{x\})= \cases{
\displaystyle\frac{1}{n}L\biggl(\tilde{V}_n^{-1}(\tau),\tilde{S}_n\biggl(x-\frac{1}{n}\biggr)\biggr),
&\quad if $ nx\in\mathbb{Z},$\cr
0, &\quad if $ nx\notin\mathbb{Z} $.
}
\]
\end{proposition}

\begin{pf}
First, we note that
\[
S_n^{-1}(S_n(x))=x+1/n \qquad  \mbox{for  all }  x \mbox{ satisfying }
nx \in\mathbb{Z}.
\]
If we use the fact that $ \{B_n(V_n^{-1}(t));t\geq0\}\}\stackrel{\mathcal{D}}{=}\{S_n(X_n(t));t\geq0\}, $ then we can\vspace*{-3pt}
see that $ \{\hat{X}_n(t);t\geq0\}\stackrel{\mathcal{D}}{=}\{
X_n(t)+1/n;t\geq0\} $.
Therefore, we see that $ \hat{X}_n $ only takes values in the lattice
$ \frac{1}{n}\mathbb{Z} $.
Moreover, we have that $ \tilde{S}_n $ and $ \tilde{V_n} $ have the
same joint distribution as
$ S_n $ and $ V_n $.
Therefore, $ \hat{X}_n=S_n^{-1}(B_n(V_n^{-1}(\cdot))) $
has the same distribution as
$ \tilde{X}_n=\tilde{S}_n^{-1}(B(\tilde{V}_n^{-1}(\cdot))) $.
From this, it also follows that $ \tilde{X}_n $ stays for all time in
the countable state space
$ \{x\in\mathbb{R};nx\in\mathbb{Z}\} $.
This implies that $ \tilde{\Gamma}_n(\tau,\{x\})=0 $ for $ nx\notin
\mathbb{Z} $.
This proves one part of the statement.

For the proof of the other part of the statement, we will need the
derivative of the function
\[
\tilde{M}(\sigma):= \frac{1}{n}L\bigl(\tilde{V}_n^{-1}(\sigma),\tilde
{S}_n(x-1/n)\bigr) .
\]
We first collect some useful facts which help to compute the derivative
of $ \tilde{M} $.\vspace*{1pt}

Since $ \tilde{S}_n $ is constant on the intervals $ [\frac
{k}{n},\frac{k+1}{n}) $ for all $ k\in\mathbb{Z} $,
we have
\begin{equation} \label{SumGleich}
\tilde{V}_n(t)=\int_\mathbb{R}L(t,\tilde{S}_n(x))\,\mathrm{d}x=\frac
{1}{n}\sum_{i\in\mathbb{Z}}L\bigl(t,\tilde{S}_n(i/n)\bigr) .
\end{equation}
Since the $ (t,x)\mapsto L(t,x) $ is jointly continuous and
non-decreasing $ \mathbb{P} $-almost surely
(see Boylan (\citeyear{Boy1964}) or Getoor and Kesten (\citeyear{GetKes1972})), it follows that $
t\mapsto\tilde{V}_n(t) $ is continuous and
non-decreasing $ \mathbb{P} \times\tilde{\mathbb{P} } $-almost
surely. This then
gives rise to
\begin{equation} \label{VnId}
\tilde{V}_n\circ\tilde{V}_n^{-1}=\operatorname{id}_{\mathbb{R}^+}
\qquad \mathbb{P} \times\tilde{\mathbb{P} }\mbox{-almost surely}.
\end{equation}
By construction, one has for all $ b\in\{\tilde{S}_n(x);x\in\mathbb
{R}\} $ that $ \tilde{S}_n^{-1}(b)=x $ is
equivalent to $ b=\tilde{S}_n(x-\frac{1}{n}) $. Moreover, one has that
$ B(\tilde{V}_n^{-1}(\sigma))\in\{\tilde{S}_n(x);x\in\mathbb{R}\}
$ for all $ \sigma\geq0 $ almost surely with
respect to $ \mathbb{P} \times\tilde{\mathbb{P} } $.
Hence,
\begin{equation} \label{equival}
\tilde{X}_n(\sigma)=\tilde{S}_n^{-1}(B(\tilde{V}_n^{-1}(\sigma
)))=x \mbox{ is equivalent to }
B(\tilde{V}_n^{-1}(\sigma))=\tilde{S}_n\biggl(x-\frac{1}{n}\biggr) .
\end{equation}
Moreover, the random variables $ \{\lambda_i^{-1};i\in\mathbb{N}\} $
are positive $ \mathbb{P} $-almost surely and
therefore
\begin{eqnarray} \label{inject}
\mbox{the restriction of }  x\mapsto\tilde{S}_n(x)  \mbox{ to the
set } \frac{1}{n}\mathbb{Z}
\mbox{ is injective almost surely with respect to } \tilde{\mathbb
{P} } .
\end{eqnarray}
Since, conditioned on $ \mathcal{A}=\sigma\{\lambda_j;j\in\mathbb{N}\}
$, the process $ X $ is a Markov process, it
follows that for
$ nx\in\mathbb{Z} $, there exist non-negative random variables $
a_1<b_1<a_2<b_2<\cdots$ with the property
\begin{eqnarray*}
\{\sigma\geq0;\tilde{X}_n(\sigma)=x \}=\bigcup_{i\in\mathbb
{N}}[a_i,b_i) \qquad
\mathbb{P} \times\tilde{\mathbb{P} }\mbox{-a.s.}
\end{eqnarray*}
This implies that for all $ \sigma_0\notin\{a_i;i\in\mathbb{N}\} $,
there exists a neighborhood
$ \mathcal{U}(\sigma_0)$ containing $ \sigma_0 $ with the property that
$ \sigma\mapsto\tilde{X}_n(\sigma)=\tilde{S}_n^{-1}(B(\tilde
{V}_n^{-1}(\sigma))) $ is constant on
$ \mathcal{U}(\sigma_0) $.
Equations (\ref{equival}) and (\ref{inject}) then imply that $ \sigma\mapsto
B(\tilde{V}_n^{-1}(\sigma)) $ must be
constant on $ \mathcal{U}(\sigma_0) $.

Therefore, for $ \sigma_0\notin\{a_i;i\in\mathbb{N}\} $ and
$ B(\tilde{V}_n^{-1}(\sigma_0))\neq\tilde{S}_n(x-\frac{1}{n}) $,
we have
$ B(\tilde{V}_n^{-1}(\sigma))\neq\tilde{S}_n(x-\frac{1}{n}) $ for
all $ \sigma$ in a neighborhood of
$ \sigma_0 $. Hence
\[
\sigma\mapsto L\bigl(\tilde{V}_n^{-1}(\sigma),\tilde{S}_n(x-1/n)\bigr)
\]
is constant in a neighborhood of $ \sigma_0$.
The previous argument and the fact that $ \tilde{X}_n $\vspace*{1pt} only jumps to
nearest neighbors in
$ \frac{1}{n}\mathbb{Z} $ leads to the fact that $ \sigma_0\notin\{
a_i;i\in\mathbb{N}\} $ and
$ B(\tilde{V}_n^{-1}(\sigma_0))=\tilde{S}_n(x-\frac{1}{n}) $ imply
the existence of a suitable $ c_0>0 $ with
the property
\[
\sigma\mapsto\frac{1}{n}\sum_{z\neq nx-1}L\bigl(\tilde{V}_n^{-1}(\sigma
),\tilde{S}_n(z/n)\bigr)=c_0
\]
in a neighborhood of $ \sigma_0$.
Therefore, we can use (\ref{VnId}) to see that $ B(\tilde
{V}_n^{-1}(\sigma_0))=\tilde{S}_n(x-\frac{1}{n}) $ implies that
\[
\sigma\mapsto\frac{1}{n}L\bigl(\tilde{V}_n^{-1}(\sigma),\tilde
{S}_n(x-1/n)\bigr)=\tilde{V}_n(\tilde{V}_n^{-1}(\sigma))-c_0=\sigma-c_0
\]
in a neighborhood of $  \sigma_0 $.
Consequently, the function
\[
\tilde{M}(\sigma):= \frac{1}{n}L\bigl(\tilde{V}_n^{-1}(\sigma),\tilde
{S}_n(x-1/n)\bigr)
\]
is differentiable for all $ \sigma\notin\{a_i;i\in\mathbb{N}\} $,
and for $ nx\in\mathbb{Z} $, we have
\[
\tilde{M}'(\sigma)= \cases{
1, &\quad if $\displaystyle B(\tilde{V}_n^{-1}(\sigma))=\tilde{S}_n\biggl(x-\frac
{1}{n}\biggr)$,\cr
0, &\quad if $\displaystyle B(\tilde{V}_n^{-1}(\sigma))\neq\tilde{S}_n\biggl(x-\frac{1}{n}\biggr)$.
}
\]
Moreover, it is possible to prove that the function $ \tilde{M} $ is
Lipschitz continuous with Lipschitz constant one.
From those properties, it follows that
\[
\int_0^\tau\mathbh{1}_{\{x\}}(\tilde{X}_n(\sigma
))\,\mathrm{d}\sigma=
\int_0^\tau\mathbh{1}_{\{\tilde{S}_n(x-1/n)\}
}(B(\tilde{V}_n^{-1}(\sigma)))\,\mathrm{d}\sigma=
\int_0^\tau\tilde{M}'(\sigma)\,\mathrm{d}\sigma=\tilde{M}(\tau) .
\]
\upqed
\end{pf}

\subsection{The convergence of the occupation times}

In this section, we investigate whether the occupation times of $
\tilde{X}_n $ converge toward the local
time of $ \tilde{X}_\ast$ in an appropriate way as $ n\rightarrow
\infty$. For this, we first need some
auxiliary results.

\begin{lemma} \label{TimeChangeConvLem}
One has $ \mathbb{P} \times\tilde{\mathbb{P} } $-almost surely that
$ \tilde{V}_n(t)
$ converges toward $ \tilde{V}_\ast(t) $
for all $ t\in\mathbb{R} $.
\end{lemma}

\begin{pf}
We fix a $ T>0 $ and define $ w_o:=\sup\{x\dvtx L(T,x)>0\} $ and $
w_u:=\inf\{x\dvtx L(T,x)>0\} $.
Those two random variables are independent of $ \tilde{\mathbb{P} } $.
We know that $ \{\tilde{S}_n(x);x\in\mathbb{R}\} $ converges toward
$ \{\tilde{W}(x);x\in\mathbb{R}\} $ with
respect to the $ J_1 $-topology $ \tilde{\mathcal{F}} $-almost surely.
We note that the local time of Brownian motion $ (x,t)\mapsto L(t,x) $
is jointly continuous $ \mathbb{P} $-almost surely
(see Boylan (\citeyear{Boy1964}) or Getoor and Kesten (\citeyear{GetKes1972})).

It follows that $ \mathbb{P} \times\tilde{\mathbb{P} } $-almost
surely $ \{L(t,\tilde
{S}_n(x));x\in\mathbb{R}\} $ converges
toward $ \{L(t,\tilde{W}(x));x\in\mathbb{R}\} $ with respect to the
$ J_1 $-topology for all
$ t\in[0,T] $.

We fix a pair $ (\omega,\tilde{\omega})\in\Omega\times\tilde
{\Omega} $ with the property that
$ \{L(t,\tilde{S}_n(x))(\omega,\tilde{\omega});x\in\mathbb{R}\} $
converges toward
$ \{L(t,\tilde{W}(x))(\omega,\tilde{\omega});x\in\mathbb{R}\} $
with respect to the $ J_1 $-topology for all
$ t\in[0,T] $.

There then exist suitable $ x_u,x_o\in\mathbb{R} $ with $ \tilde
{W}(x_u)\leq w_u $ and $ \tilde{W}(x_o)\geq w_o $,
and there exists a sequence of increasing, absolutely continuous,
surjective Lipschitz maps
$ \lambda_n\dvtx [x_u,x_o]\rightarrow[x_u,x_o] $ with the properties
\[
\sup_{x\in[x_u,x_o]} |L(t,\tilde{W}(x))-L(t,\tilde{S}_n(\lambda
_n(x))) |\longrightarrow0
\qquad \mbox{as }  n\rightarrow\infty
\]
and
\[
\operatorname{esssup}\limits_{x\in[x_u,x_o]} |\lambda_n'(x)-1 |\longrightarrow0
\qquad \mbox{as }  n\rightarrow\infty.
\]
We should emphasize that the derivative of the function $ \lambda_n $
may not exist everywhere.
However, those points where it does not exist form a zero set since $
\lambda_n $ is an absolutely
continuous Lipschitz function.

By a change of variables for all $ t\in[0,T] $, one then has
\begin{eqnarray*}
&& \int_{x_u}^{x_o} L(t,\tilde{S}_n(x))\,\mathrm{d}x-\int_{x_u}^{x_o}
L(t,\tilde{S}_n(\lambda_n(x)))\,\mathrm{d}x \\
&&\quad= \int_{x_u}^{x_o} L(t,\tilde{S}_n(x)) \biggl(1-\frac{1}{\lambda
_n'(\lambda_n^{-1}(x))} \biggr)\,\mathrm{d}x
+\mathrm{O} \Bigl(\sup_{x\in[x_u,x_o]}|\lambda_n(x)-x| \Bigr).
\end{eqnarray*}
It follows from the assumptions on the sequence $ \lambda_n $ that the
above difference converges toward zero.
Further, for all $ t\in[0,T] $, we have that
\[
\int_{\mathbb{R}} L(t,\tilde{S}_n(\lambda_n(x)))\,\mathrm{d}x\longrightarrow
\int_{\mathbb{R}} L(t,\tilde{W}(x))\,\mathrm{d}x\qquad \mbox{as } n\rightarrow
\infty.
\]
Hence, one has $ \mathbb{P} \times\tilde{\mathbb{P} } $-almost
surely that $ \tilde
{V}_n(t) $ converges toward
$ \tilde{V}_\ast(t) $ for all $ t\in[0,T] $. Thus, for every $ T>0
$, we obtain an zero set $ N_T $ in
$ \Omega\times\tilde{\Omega} $ where this convergence does not
hold. The lemma now follows since the union
\[
N_\infty:=\bigcup_{T\in\mathbb{N}}U_T
\]
is also a zero set with respect to $ \mathbb{P} \times\tilde{\mathbb
{P} } $.
\end{pf}

Let $ f\dvtx\mathbb{R}\rightarrow\mathbb{R} $ be a function. We call $
\tau\in f(\mathbb{R}) $ a
\emph{critical value} for $ f $ if there exist at least two distinct
points $ t_1,t_2\in\mathbb{R} $ such that
$ f(t_1)=f(t_2)=\tau$.
Further, we call a point $ \tau\in f(\mathbb{R}) $ a \emph{regular
value} for $ f $ if it is not a critical
value.
It is straightforward to see that the preimages of critical values
contain an open interval
if the function $ f $ is non-decreasing.
This implies that the set of critical values of a non-decreasing
function is at most countable.

\begin{lemma}\label{InversTimeChangeConvLem}
One has $ \mathbb{P} \times\tilde{\mathbb{P} } $-almost surely that
$ \tilde
{V}_n^{-1}(\tau) $ converges toward
$ \tilde{V}_\ast^{-1}(\tau) $ for all regular values $ \tau$ of $
\tilde{V}_\ast$.
\end{lemma}

\begin{pf}
We note that $ \mathbb{P} $-almost surely the local time $ L(t,x) $ of the
Brownian motion $ B $ is continuous and
non-decreasing in $ t $ for all $ x\in\mathbb{R} $ (see Boylan (\citeyear{Boy1964})
or Getoor and Kesten (\citeyear{GetKes1972})) for
the continuity).
It follows that $ \mathbb{P} \times\tilde{\mathbb{P} } $-almost
surely the function
\[
t\mapsto\tilde{V}_\ast(t):=\int_\mathbb{R} L(t,x)m_\ast(\mathrm{d}x)
\]
is continuous and non-decreasing.

Therefore, $ \mathbb{P} \times\tilde{\mathbb{P} } $-almost surely
the function
$ \tilde{V}_\ast^{-1}(\tau):=\inf\{t;\tilde{V}(t)>\tau\} $ is
strictly increasing and right-continuous.

We use Lemma \ref{TimeChangeConvLem} to fix a pair $ (\omega,\tilde
{\omega})\in\Omega\times\tilde{\Omega} $
with the properties that:
\begin{longlist}[(ii)]
\item[(i)] $ \tau\mapsto\tilde{V}_\ast^{-1}(\tau) $ is strictly
increasing and right-continuous;
\item[(ii)] $ \tilde{V}_n(t) $ converges toward $ \tilde{V}_\ast(t) $ for
all $ t\geq0 $.
\end{longlist}
Since the set where $ \tilde{V}_\ast$ is not continuous is countable,
the set where $ \tilde{V}_\ast$
is continuous is dense in $ [0,\infty) $.

We denote by $ K $ the set of critical values of $ \tilde{V}_\ast$.
As was pointed out before, $ K $ is at most countable.
For an arbitrary point $ \tau\in[0,\infty)\cap K^c $
and for any $ \epsilon>0 $,
one can find points
$ t_{\epsilon,0}, t_{\epsilon,1}\in(\tilde{V}_\ast^{-1}(\tau
)-\epsilon,\tilde{V}_\ast^{-1}(\tau)) $ and
$ t_{\epsilon,2}, t_{\epsilon,3}\in(\tilde{V}_\ast^{-1}(\tau
),\tilde{V}_\ast^{-1}(\tau)+\epsilon) $
with the property
\[
\tilde{V}_\ast(t_{\epsilon,0})<\tilde{V}_\ast(t_{\epsilon
,1})<\tau<\tilde{V}_\ast(t_{\epsilon,2})
<\tilde{V}_\ast(t_{\epsilon,3}) .
\]
We can now choose a $ \delta>0$ such that
\[
\tilde{V}_\ast(t_{\epsilon,0})+\delta<\tilde{V}_\ast(t_{\epsilon
,1})-\delta
<\tilde{V}_\ast(t_{\epsilon,1})+\delta<\tau<\tilde{V}_\ast
(t_{\epsilon,2})-\delta
<\tilde{V}_\ast(t_{\epsilon,2})+\delta
<\tilde{V}_\ast(t_{\epsilon,3})-\delta.
\]
Since $ \tilde{V}_n $ converges toward $ \tilde{V}_\ast$ in all
points where $ \tilde{V}_\ast$ is continuous,
there exists an $ n_0\in\mathbb{N} $ such that for all $ n\geq n_0 $,
we have
\[
\tilde{V}_n(t_{\epsilon,0})<\tilde{V}_\ast(t_{\epsilon,0})+\delta
<\tilde{V}_\ast(t_{\epsilon,1})-\delta
<\tilde{V}_n(t_{\epsilon,1})<\tilde{V}_\ast(t_{\epsilon,1})+\delta
<\tau
\]
and
\[
\tau<\tilde{V}_\ast(t_{\epsilon,2})-\delta<\tilde
{V}_n(t_{\epsilon,2})
<\tilde{V}_\ast(t_{\epsilon,2})+\delta
<\tilde{V}_\ast(t_{\epsilon,3})-\delta<\tilde{V}_n(t_{\epsilon
,3}) .
\]
By definition of $ t_{\epsilon,0}$, we have that $ z\leq\tilde
{V}_\ast^{-1}(\tau)-\epsilon$ implies
$ z\leq t_{\epsilon,0} $.
From monotonicity and the first of both inequalities above, it follows that
\[
\tilde{V}_n(z)\leq\tilde{V}_n(t_{\epsilon,0})\leq\tilde{V}_\ast
(t_{\epsilon,0})
+\delta<\tilde{V}_\ast(t_{\epsilon,1}) .
\]
We have thus seen that $ z\leq\tilde{V}_\ast^{-1}(\tau)-\epsilon$ implies
$ \tilde{V}_n(z)<\tilde{V}_\ast(t_{\epsilon,1}) $.
If we reverse the implication, then we obtain that $ \tilde
{V}_n(z)\geq\tilde{V}_\ast(t_{\epsilon,1}) $ implies
$ z>\tilde{V}_\ast^{-1}(\tau)-\epsilon$. From this implication, it
follows that
\[
\tilde{V}_n^{-1}(\tilde{V}_\ast(t_{\epsilon,1}))=\inf\{z\dvtx\tilde
{V}_n(z)>\tilde{V}_\ast(t_{\epsilon,1})\}
>\tilde{V}_\ast^{-1}(\tau)-\epsilon.
\]
For $ z=t_{\epsilon,3} $, we have $ \tilde{V}_n(z)=\tilde
{V}_n(t_{\epsilon,3})>\tilde{V}_\ast(t_{\epsilon,2}) $.
In other words, there exists a $ z<\tilde{V}_\ast^{-1}(\tau
)+\epsilon$ with
$ \tilde{V}_n(z)>\tilde{V}_\ast(t_{\epsilon,2}) $.
This proves that
\[
\tilde{V}_\ast^{-1}(\tau)+\epsilon>\tilde{V}_n^{-1}(\tilde
{V}_\ast(t_{\epsilon,2})) .
\]
Altogether, we have proven that for all $ n\geq n_0 $,
\[
\tilde{V}_\ast^{-1}(\tau)-\epsilon<\tilde{V}_n^{-1}(\tilde
{V}_\ast(t_{\epsilon,1}))
<\tilde{V}_n^{-1}(\tilde{V}_\ast(t_{\epsilon,2}))<\tilde{V}_\ast
^{-1}(\tau)+\epsilon.
\]
By monotonicity, for all $ n\geq n_0 $ and all
$ \tau'\in[\tilde{V}_\ast(t_{\epsilon,1}),\tilde{V}_\ast
(t_{\epsilon,2})], $ one has
\[
\tilde{V}_\ast^{-1}(\tau)-\epsilon<\tilde{V}_n^{-1}(\tau')<\tilde
{V}_\ast^{-1}(\tau)+\epsilon.
\]
Since $ \tau\in[\tilde{V}_\ast(t_{\epsilon,1}),\tilde{V}_\ast
(t_{\epsilon,2})] $, the proof is complete.
\end{pf}

\begin{lemma} \label{RegValLem1}
For all $ \tau\geq0 $, one has that $ \tau$ is a regular value of $
\tilde{V}_\ast$ almost surely with
respect to $ \mathbb{P} \times\tilde{\mathbb{P} } $.
\end{lemma}

\begin{pf}
By the invariance properties of Brownian motion, we have that for all $
\gamma>0 $,
\[
\{L(t,w);w\in\mathbb{R},t\geq0\}\stackrel{\mathcal{D}}{=}
\{\gamma^{-1}L(\gamma^2t,\gamma w);w\in\mathbb{R},t\geq0\} .
\]
By the invariance of the $ \alpha$-stable L\'{e}vy process, we have that
\begin{eqnarray*}
\{L(t,\tilde{W}(x));x\in\mathbb{R},t\geq0\}&\stackrel{\mathcal{D}}{=}&
\{\gamma^{-1}L(\gamma^2t,\gamma\tilde{W}(x));x\in\mathbb{R},t\geq
0\} \\
& \stackrel{\mathcal{D}}{=}& \{\gamma^{-1}L(\gamma^2t,\tilde{W}(\gamma
^\alpha x));x\in\mathbb{R},t\geq0\} .
\end{eqnarray*}
Substitution then yields
\begin{eqnarray*}
\biggl\{\int_{\mathbb{R}}L(t,\tilde{W}(x))\,\mathrm{d}x;t\geq0 \biggr\}
& \stackrel{\mathcal{D}}{=}& \biggl\{\gamma^{-1}\int_{\mathbb{R}}L(\gamma
^2t,\tilde{W}(\gamma^\alpha x))\,\mathrm{d}x;t\geq0 \biggr\} \\
& \stackrel{\mathcal{D}}{=}& \biggl\{\gamma^{-1-\alpha}\int_{\mathbb
{R}}L(\gamma^2t,\tilde{W}(x))\,\mathrm{d}x;t\geq0 \biggr\} .
\end{eqnarray*}
By definition, this means that
\[
\{\tilde{V}_\ast(t);t\geq0\} \stackrel{\mathcal{D}}{=}
\{\gamma^{-1-\alpha}\tilde{V}_\ast(\gamma^2t);t\geq0\} .
\]
We define $ \ell_\ast$ to be the image measure of the Lebesgue
measure $ \ell$ with respect
$ \tilde{V}_\ast$.
The previous considerations imply that
\[
\ell_\ast(\mathrm{d}t)\stackrel{\mathcal{D}}{=}\gamma^2\ell_\ast(\gamma^{-1-\alpha}\,\mathrm{d}t) .
\]
This identity implies that no $ \tau>0 $ satisfies $ \ell_\ast(\{
\tau\})>0 $ with a positive probability
with respect to $ \mathbb{P} \times\tilde{\mathbb{P} } $.
To a critical value $ \tau$ corresponds an interval where $ t\mapsto
\tilde{V}_\ast$ is constant, which
implies that $ \ell_\ast(\{\tau\})>0 $. For a particular point $
\tau>0 $, this cannot happen with positive
probability. This finishes the proof of the statement.
\end{pf}

\begin{proposition} \label{OkTimeLokTimeConvProp}
For all $ \tau\geq0 $, the sequence of functions $ x\mapsto L(\tilde
{V}_n^{-1}(\tau),\tilde{S}_n(x+1/n)) $
converges toward the function $ x\mapsto L(\tilde{V}_\ast^{-1}(\tau
),\tilde{W}(x)) $ in the $ J_1 $-topology
$ \mathbb{P} \times\tilde{\mathbb{P} } $-almost surely.
\end{proposition}

\begin{pf}
It is known that $ \tilde{S}_n $ converges toward $ \tilde{W} $ in
the $ J_1 $-topology almost surely with
respect to $ \tilde{\mathbb{P} } $.
Moreover, by Lemmas \ref{InversTimeChangeConvLem} and \ref{RegValLem1}, for all $ \tau\geq0 $, the
sequence $ \tilde{V}_n^{-1}(\tau) $ converges toward $ \tilde
{V}_\ast^{-1}(\tau) $ almost surely
with respect to $ \mathbb{P} \times\tilde{\mathbb{P} } $. The
proposition follows
since it is well known that
$ (t,x)\mapsto L(t,x) $ is jointly continuous $ \mathbb{P} $-almost
surely; see
Boylan (\citeyear{Boy1964}) or
Getoor and Kesten (\citeyear{GetKes1972}).
\end{pf}

\begin{lemma} \label{LokalNullMenge}
For all $ k\in\mathbb{N} $, $ \theta_1,\ldots,\theta_k\in\mathbb{R}
$ and all $ \tau_1,\ldots,\tau_k\geq0 $, the set
\[
\mathcal{C}:=
\Biggl\{c>0\dvtx\ell\Biggl(x\in\mathbb{R}; \Biggl|\sum_{i=1}^k\theta_iL(\tilde{V}_\ast
^{-1}(\tau_i),\tilde{W}(x)) \Biggr|=c \Biggr)>0 \Biggr\}
\]
is countable $ \mathbb{P} \times\tilde{\mathbb{P} } $-almost
surely, where $ \ell$
denotes the Lebesgue measure on
$ \mathbb{R} $.
\end{lemma}
\begin{pf}
It is well known that $ x\mapsto\tilde{W}(x) $ is strictly increasing
$ \tilde{\mathbb{P} } $-almost surely.
For $ c>0 $, we define the level-sets
\[
\mathcal{N}_c:= \Biggl\{w\in\mathbb{R}; \Biggl|\sum_{i=1}^k\theta_i
L(\tilde{V}_\ast^{-1}(\tau_i),w) \Biggr|=c \Biggr\} .
\]
Fix a strictly increasing path $ f\dvtx x\mapsto\tilde{W}(x) $ and assume
that there exist an uncountable number of
$ c>0 $ with the property that $ \ell(f^{-1}(\mathcal{N}_c))>0 $.
For $ c\neq c' $, the sets $ f^{-1}(\mathcal{N}_c) $ and $ f^{-1}(\mathcal{N}_{c'}) $ are disjoint.
We would obtain an uncountable number of disjoint sets with positive
Lebesgue measure.
This is, of course, not possible.
\end{pf}

\begin{proposition} \label{CardTowardMeasureProp}
For all $ k\in\mathbb{N} $, $ \theta_1,\ldots,\theta_k\in\mathbb
{R} $ and all $ \tau_1,\ldots,\tau_k\geq0 $,
one has $ \mathbb{P} \times\tilde{\mathbb{P} } $-almost surely that
\begin{eqnarray*}
&&\frac{1}{n}\operatorname{card}
\Biggl\{x\in\mathbb{Z}\dvtx n \Biggl|\sum_{i=1}^k\theta_i\tilde{\Gamma}_n(\tau
_i,\{x/n\}) \Biggr|>c \Biggr\}\\
&&\quad\longrightarrow
\ell\Biggl(x\in\mathbb{R}\dvtx \Biggl|\sum_{i=1}^k\theta_i\tilde{L}_\ast(\tau
_i,x) \Biggr|>c \Biggr) \qquad
\mbox{as } n\rightarrow\infty
\end{eqnarray*}
for all but a countable number of $ c>0 $.
\end{proposition}
\begin{pf}
We can find a $ K>0 $ such that $ \{y\in\mathbb{R}\dvtx L(\tau_i,y)\neq0
\ {\rm for\ all}\ i=1,\ldots,k \} $ is a subset of the
interval $ (\tilde{W}(-K),\tilde{W}(K)) $. By Propositions \ref
{OkTimePro}, \ref{OkTimeLokTimeConvProp} and
Corollary \ref{LokTimeKor1}, the
sequence
\begin{eqnarray*}
\tilde{A}_n(x)&:=&n \Biggl|\sum_{i=1}^k\theta_i\tilde{\Gamma}_n(\tau_i,\{
x\}) \Biggr|\\
&\hspace*{3pt}=& \Biggl|\sum_{i=1}^k\theta_iL\bigl(\tilde{V}_n^{-1}(\tau_i),\tilde
{S}_n(x-1/n)\bigr) \Biggr|
\end{eqnarray*}
converges $ \mathbb{P} \times\tilde{\mathbb{P} } $-almost surely in
the $ J_1$-topology toward
\[
\tilde{A}_\ast(x):= \Biggl|\sum_{i=1}^k\theta_i\tilde{L}_\ast(\tau
_i,x) \Biggr|
= \Biggl|\sum_{i=1}^k\theta_iL(\tilde{V}_\ast^{-1}(\tau_i),\tilde
{W}(x)) \Biggr| .
\]
There then exists a sequence of continuous increasing maps $ \lambda
_n\dvtx[-K,K]\rightarrow[-K,K] $ such that
\[
\sup_{x\in[-K,K]} |\tilde{A}_\ast(x)-\tilde{A}_n\circ\lambda
_n(x) |\longrightarrow0
\qquad \mbox{as }  n\rightarrow\infty
\]
and such that each $ \lambda_n $ is Lipschitz continuous and satisfies
\[
\operatorname{esssup}\limits_{x\in[-K,K]} |\lambda_n'(x)-1 |\longrightarrow0 .
\]
We should emphasize that the derivative of the function $ \lambda_n $
may not exist everywhere.
However, those points where the derivative does not exist form a zero
set since $ \lambda_n $ is an
absolutely continuous Lipschitz function. We note that for suitably
large $ n\in\mathbb{N}$, one has
\begin{eqnarray*}
&&\frac{1}{n}\operatorname{card} \Biggl\{x\in\mathbb{R};
\Biggl|\sum_{i=1}^k\theta_iL\bigl(\tilde{V}_n^{-1}(\tau_i),\tilde
{S}_n(x-1/n)\bigr) \Biggr|>c \Biggr\}\\
&&\quad=\ell\bigl(x\in[-K,K];\tilde{A}_n(x)>c \bigr)=\int_{-K}^K\mathbh{1}_{(c,\infty)}(\tilde{A}_n(x))\,\mathrm{d}x .
\end{eqnarray*}
It then follows that
\begin{eqnarray*}
&& \frac{1}{n}\operatorname{card} \Biggl\{x\in[-K,K];
n \Biggl|\sum_{i=1}^k\theta_i\tilde{\Gamma}_n(\tau_i,\{x\}) \Biggr|>c \Biggr\}
-\int_{-K}^K\mathbh{1}_{(c,\infty)}(\tilde{A}_n(\lambda
_n(x)))\,\mathrm{d}x\\
&&\quad= \int_{-K}^K\mathbh{1}_{(c,-\infty)}(\tilde{A}_n(x))\,\mathrm{d}x
\biggl(1-\frac{1}{\lambda_n'(\lambda_n^{-1}(x))} \biggr)\,\mathrm{d}x
+\mathrm{O} \Bigl(\sup_{x\in[-K,K]}|\lambda_n(x)-x| \Bigr).
\end{eqnarray*}
By the assumptions on the sequence $ \{\lambda_n;n\in\mathbb{N}\} $,
the previous
difference converges toward zero. Furthermore,
\[
\int_{-K}^K\mathbh{1}_{(c,\infty)}(\tilde{A}_n(\lambda
_n(x)))\,\mathrm{d}x\longrightarrow
\int_{-K}^K\mathbh{1}_{(c,\infty)}(\tilde{A}_\ast
(x))\,\mathrm{d}x\qquad \mbox{as } n\rightarrow\infty
\]
whenever the set $ \{x\in[-K,K];\tilde{A}_\ast(s)=c\} $ is a zero
set with respect to the Lebesgue measure
$ \ell$ on $ \mathbb{R} $. Since this was proven in Lemma \ref
{LokalNullMenge}, the statement of the
proposition follows.
\end{pf}

Subsequently, we will make use of the following notation:
\[
A_n^+:= \Biggl\{x\in\mathbb{Z}\dvtx\sum_{i=1}^k\theta_i\tilde{\Gamma
}_n(\tau_i,\{x/n\})>0 \Biggr\} ,\qquad
A_n^-:= \Biggl\{x\in\mathbb{Z}\dvtx\sum_{i=1}^k\theta_i\tilde{\Gamma
}_n(\tau_i,\{x/n\})<0 \Biggr\}
\]
and
\[
A^+:= \Biggl\{x\in\mathbb{R}\dvtx\sum_{i=1}^k\theta_i\tilde{L}_\ast(\tau
_i,x)>0 \Biggr\} ,\qquad
A^-:= \Biggl\{x\in\mathbb{R}\dvtx\sum_{i=1}^k\theta_i\tilde{L}_\ast(\tau
_i,x)<0 \Biggr\} .
\]
Later, we will need the following version of Proposition \ref
{CardTowardMeasureProp}.

\begin{proposition} \label{SignedCardTowardMeasureProp}
For all $ k\in\mathbb{N} $, $ \theta_1,\ldots,\theta_k\in\mathbb
{R} $ and all $ \tau_1,\ldots,\tau_k\geq0 $, one has
$ \mathbb{P} \times\tilde{\mathbb{P} } $-almost surely that
\[
\frac{1}{n}\operatorname{card}
\Biggl\{x\in\mathbb{Z}\cap A_n^\pm\dvtx
n \Biggl|\sum_{i=1}^k\theta_i\tilde{\Gamma}_n(\tau_i,\{x/n\}) \Biggr|>c \Biggr\}
\longrightarrow
\ell\Biggl(x\in\mathbb{R}\cap A^\pm\dvtx \Biggl|\sum_{i=1}^k\theta_i\tilde
{L}_\ast(\tau_i,x) \Biggr|>c \Biggr)
\]
for all but a countable number of $ c>0 $.
\end{proposition}

\begin{pf}
The proof uses essentially the same arguments as the proof of
Proposition \ref{CardTowardMeasureProp}.
\end{pf}

\begin{remark*}
With the same proof as for Proposition \ref
{CardTowardMeasureProp}, we can show that
$ \mathbb{P} \times\tilde{\mathbb{P} } $-almost surely
\[
\frac{1}{n}\operatorname{card} \bigl\{x\in\mathbb{Z}\dvtx n^2\tilde{\Gamma}_n^2(\tau
_i,\{x/n\})>c \bigr\}\longrightarrow
\ell\bigl(x\in\mathbb{R}\dvtx\tilde{L}_\ast^2(\tau_i,x)>c \bigr) \qquad
\mbox{as } n\rightarrow\infty
\]
for all but a countable number of $ c>0 $.
\end{remark*}

\subsection{A useful lemma on integrated powers of local time}

\begin{lemma} \label{PrinceLem}
For $ \tau_1,\ldots,\tau_k\geq0 $ and $ \theta_1,\ldots,\theta
_k\in\mathbb{R} $, the two sequences of random variables
\begin{eqnarray*}
&&n^{\beta-1}\sum_{x\in\mathbb{Z}} \Biggl|\sum_{i=1}^k\theta_i\tilde
{\Gamma}_n(\tau_i,\{x/n\}) \Biggr|^\beta\quad
\mbox{and }\\
&&   n^{\beta-1}\sum_{x\in\mathbb{Z}}
\Biggl( \Biggl|\sum_{i=1}^k\theta_i\tilde{\Gamma}_n(\tau_i,\{x/n\}) \Biggr|^\beta
\operatorname{sgn} \Biggl(\sum_{i=1}^k\theta_i\tilde{\Gamma}_n(\tau_i,\{x/n\}) \Biggr)
\Biggr)
\end{eqnarray*}
converge $ \mathbb{P} \times\tilde{\mathbb{P} } $-almost surely
toward the respective
random variables
\[
\int_{-\infty}^\infty\Biggl|\sum_{i=1}^k\theta_i\tilde{L}_\ast(\tau
_i,x) \Biggr|^\beta \,\mathrm{d}x \quad
\mbox{and}\quad
\int_{-\infty}^\infty\Biggl( \Biggl|\sum_{i=1}^k\theta_i\tilde{L}_\ast(\tau
_i,x) \Biggr|^\beta
\operatorname{sgn} \Biggl(\sum_{i=1}^k\theta_i\tilde{L}_\ast(\tau_i,x) \Biggr) \Biggr)\,\mathrm{d}x .
\]
\end{lemma}

\begin{pf}
We use the layer cake representation of the integrals (see Lieb and
Loss (\citeyear{LieLos2001})) to write
\[
\sum_{x\in\mathbb{Z}} \Biggl|\sum_{i=1}^k\theta_in\tilde{\Gamma
}_n(\tau_i,\{x/n\}) \Biggr|^\beta=
\beta\int_0^\infty c^{\beta-1}\operatorname{card}
\Biggl\{x\in\mathbb{Z}\dvtx n \Biggl|\sum_{i=1}^k\theta_i\tilde{\Gamma}_n(\tau
_i,\{x/n\}) \Biggr|>c \Biggr\}\,\mathrm{d}c
\]
and
\[
\int_{-\infty}^\infty\Biggl|\sum_{i=1}^k\theta_i\tilde{L}_\ast(\tau
_i,x) \Biggr|^\beta \,\mathrm{d}x
=\beta\int_0^\infty c^{\beta-1}\ell\Biggl(x\in\mathbb{R}\dvtx
\Biggl|\sum_{i=1}^k\theta_i\tilde{L}_\ast(\tau_i,x) \Biggr|>c \Biggr)\,\mathrm{d}c .
\]
We note that the convergence of $ \tilde{V}_n^{-1}(\tau_i) $ toward $
\tilde{V}_\ast^{-1}(\tau_i) $ and
the fact that $ t\mapsto L(t,y) $ is increasing for every $ y\in
\mathbb{R} $ imply that there exists an
$ n_0\in\mathbb{N} $ with
\[
L(\tilde{V}_n^{-1}(\tau_i),y)\leq
L\bigl(\tilde{V}_\ast^{-1}(\tau_i)+1,y\bigr)\qquad
\mbox{for all } y\in\mathbb{R}, 1\leq i\leq k, n\geq n_0 .
\]
Moreover, for all $ i\in\{1,\ldots,k\} $, the functions $ y\mapsto
L(\tilde{V}_\ast^{-1}(\tau_i)+1,y) $ are
continuous and their supports are contained in $ [-K,K] $ for a
suitable $ K>0 $.
Hence, there exists a $ C>0 $ such that for $ n\geq n_0 $, one has
\begin{eqnarray*}
n \Biggl|\sum_{i=1}^k\theta_i\tilde{\Gamma}_n(\tau_i,\{x/n\}) \Biggr|
&\leq&\Biggl|\sum_{i=1}^k\theta_iL\bigl(\tilde{V}_n^{-1}(\tau_i),\tilde
{S}_n\bigl((x-1)/n\bigr)\bigr) \Biggr|\\
&\leq&\sum_{i=1}^k\theta_i\sup_{y\in\mathbb{R}}L\bigl(\tilde{V}_\ast
^{-1}(\tau_i)+1,y\bigr)\leq C.
\end{eqnarray*}
This implies that all of the functions
\[
c\mapsto\operatorname{card} \Biggl\{x\in\mathbb{Z}\dvtx
n \Biggl|\sum_{i=1}^k\theta_i\tilde{\Gamma}_n(\tau_i,\{x/n\}) \Biggr|>c \Biggr\}
\]
have support contained in $ [0,C]$.
Moreover, for all $ c>0 $, we have
\begin{eqnarray*}
\operatorname{card}
\Biggl\{x\in\mathbb{Z}\dvtx n \Biggl|\sum_{i=1}^k\theta_i\tilde{\Gamma}_n(\tau
_i,\{x/n\}) \Biggr|>c \Biggr\}
\leq\operatorname{card} \bigl\{x\in\mathbb{Z}\dvtx-K\leq\tilde{S}_n\bigl((x-1)/n\bigr)\leq K \bigr\}.
\end{eqnarray*}
Since
\[
\ell\bigl(x;\tilde{W}(x)\in\{-K,K\} \bigr)=0
\]
and since $ \tilde{S}_n $ converges toward $ \tilde{W} $ with respect
to the Skorohod metric, we have that
\[
\frac{1}{n}\operatorname{card} \bigl\{x\in\mathbb{Z}\dvtx-K\leq\tilde
{S}_n\bigl((x-1)/n\bigr)\leq K \bigr\}
\longrightarrow\ell\bigl(x\in\mathbb{R}\dvtx-K\leq\tilde{W}(x)\leq K \bigr).
\]
This implies that there exists an $ R>0 $ such that for all $ n\in
\mathbb{N} $ and all $ c>0 $, we have
\[
\frac{1}{n}\operatorname{card} \Biggl\{x\in\mathbb{Z}\dvtx
n \Biggl|\sum_{i=1}^k\theta_i\tilde{\Gamma}_n(\tau_i,\{x/n\}) \Biggr|>c \Biggr\}
\leq R.
\]
The first statement of the lemma then follows from dominated
convergence and Proposition
\ref{CardTowardMeasureProp}.\\[0.5mm]
The second statement is proved in the same way by separating the
positive and the negative parts of the
integrals and using the statements from Proposition \ref
{SignedCardTowardMeasureProp} instead of
Proposition \ref{CardTowardMeasureProp}.
\end{pf}

\begin{proposition} \label{PrinceProp}
For $ \tau_1,\ldots,\tau_k\geq0 $ and $ \theta_1,\ldots,\theta
_k\in\mathbb{R} $, the two sequences of random variables
\begin{eqnarray*}
&&n^{\beta-1}\sum_{x\in\mathbb{Z}} \Biggl|\sum_{i=1}^k\theta_i\Gamma
_n(\tau_i,\{x/n\}) \Biggr|^\beta\quad
\mbox{and}\\
&&n^{\beta-1}\sum_{x\in\mathbb{Z}} \Biggl( \Biggl|\sum_{i=1}^k\theta_i\Gamma
_n(\tau_i,\{x/n\}) \Biggr|^\beta
\operatorname{sgn} \Biggl(\sum_{i=1}^k\theta_i\Gamma_n(\tau_i,\{x/n\}) \Biggr) \Biggr)
\end{eqnarray*}
converge jointly in distribution toward the respective random variables
\[
\int_{-\infty}^\infty\Biggl|\sum_{i=1}^k\theta_iL_\ast(\tau_i,x)
\Biggr|^\beta \,\mathrm{d}x \quad
\mbox{and}\quad
\int_{-\infty}^\infty\Biggl( \Biggl|\sum_{i=1}^k\theta_iL_\ast(\tau_i,x)
\Biggr|^\beta
\operatorname{sgn} \Biggl(\sum_{i=1}^k\theta_iL_\ast(\tau_i,x) \Biggr) \Biggr)\,\mathrm{d}x .
\]
\end{proposition}

\begin{pf}
We know that
\[
\{L_\ast(t,x);t\geq0,x\in\mathbb{R} \} \stackrel{\mathcal{D}}{=}
\{\tilde{L}_\ast(t,x);t\geq0,x\in\mathbb{R} \}
\]
and
\[
\{S_n^{-1}(B_n(V_n^{-1}(t)));t\geq0 \} \stackrel{\mathcal{D}}{=}
\{\tilde{S}_n^{-1}(B(\tilde{V}_n^{-1}(t)));t\geq0 \} .
\]
Therefore, by Lemma \ref{PrinceLem}, the sequences of random variables
\begin{eqnarray*}
&&n^{\beta-1}\sum_{x\in\mathbb{Z}} \Biggl|\sum_{i=1}^k\theta_i\hat
{\Gamma}_n(\tau_i,\{x/n\}) \Biggr|^\beta\quad
\mbox{and}\\
&&n^{\beta-1}\sum_{x\in\mathbb{Z}} \Biggl( \Biggl|\sum_{i=1}^k\theta_i\hat
{\Gamma}_n(\tau_i,\{x/n\}) \Biggr|^\beta
\operatorname{sgn} \Biggl(\sum_{i=1}^k\theta_i\hat{\Gamma}_n(\tau_i,\{x/n\}) \Biggr) \Biggr)
\end{eqnarray*}
converge jointly in distribution toward the respective random variables
\[
\int_{-\infty}^\infty\Biggl|\sum_{i=1}^k\theta_iL_\ast(\tau_i,x)
\Biggr|^\beta \,\mathrm{d}x \quad
\mbox{and}\quad
\int_{-\infty}^\infty\Biggl( \Biggl|\sum_{i=1}^k\theta_iL_\ast(\tau_i,x)
\Biggr|^\beta
\operatorname{sgn} \Biggl(\sum_{i=1}^k\theta_iL_\ast(\tau_i,x) \Biggr) \Biggr)\,\mathrm{d}x .
\]
Moreover, $ S_n^{-1}(S_n(x/n))=(x+1)/n $ for all $ x\in\mathbb{Z} $.
This implies that
\[
\hat{X}_n(\tau)\stackrel{\mathcal{D}}{=}S_n^{-1}(S_n(X_n(\tau
)))=X_n(\tau)+1/n .
\]
Hence, we have $ \hat{\Gamma}_n(\tau,\{x/n\})\stackrel{\mathcal{D}}{=}\Gamma_n(\tau,\{(x+1)/n\}) $ for all $ x\in\mathbb{Z} $.
Therefore,
\[
n^{\beta-1}\sum_{x\in\mathbb{Z}} \Biggl|\sum_{i=1}^k\theta_i\hat
{\Gamma}_n(\tau_i,\{x/n\}) \Biggr|^\beta\stackrel{\mathcal{D}}{=}
n^{\beta-1}\sum_{x\in\mathbb{Z}} \Biggl|\sum_{i=1}^k\theta_i\Gamma
_n(\tau_i,\{x/n\}) \Biggr|^\beta
\]
and
\begin{eqnarray*}
&& n^{\beta-1}\sum_{x\in\mathbb{Z}} \Biggl( \Biggl|\sum_{i=1}^k\theta_i\hat
{\Gamma}_n(\tau_i,\{x/n\}) \Biggr|^\beta
\operatorname{sgn} \Biggl(\sum_{i=1}^k\theta_i\hat{\Gamma}_n(\tau_i,\{x/n\}) \Biggr) \Biggr)
\\
&&\quad\stackrel{\mathcal{D}}{=} n^{\beta-1}\sum_{x\in\mathbb{Z}} \Biggl( \Biggl|\sum
_{i=1}^k\theta_i\Gamma_n(\tau_i,\{x/n\}) \Biggr|^\beta
\operatorname{sgn} \Biggl(\sum_{i=1}^k\theta_i\Gamma_n(\tau_i,\{x/n\}) \Biggr) \Biggr) .
\end{eqnarray*}
This proves the proposition.
\end{pf}

For the sequel, we define the occupation time
\[
\Gamma(t,A):=\int_0^t\mathbh{1}_{A}(X(s))\,\mathrm{d}s
\]
of the process $ X $ in the measurable set $ A\subset\mathbb{R} $.
Consequently, we have
\[
\Xi(t)=\sum_x\Gamma(t,\{x\})\xi(x) .
\]
We will use this fact and the following corollary in the proofs of the
next section.
\begin{corollary} \label{PrinceKor}
For $ \tau_1,\ldots,\tau_k\geq0 $ and $ \theta_1,\ldots,\theta
_k\in\mathbb{R} $, the two sequences of random variables
\begin{eqnarray*}
&&n^{-1-{\beta/\alpha}}\sum_{x\in\mathbb{Z}} \Biggl|\sum
_{i=1}^k\theta_i\Gamma(k_n\tau_i,\{x\}) \Biggr|^\beta\quad
\mbox{and}\\
&&n^{-1-{\beta/\alpha}}\sum_{x\in\mathbb{Z}} \Biggl( \Biggl|\sum
_{i=1}^k\theta_i\Gamma(k_n\tau_i,\{x\}) \Biggr|^\beta
\operatorname{sgn} \Biggl(\sum_{i=1}^k\theta_i\Gamma(k_n\tau_i,\{x\}) \Biggr) \Biggr)
\end{eqnarray*}
converge jointly in distribution toward the respective random variables
\begin{eqnarray*}
&&\int_{-\infty}^\infty\Biggl|\sum_{i=1}^k\theta_iL_\ast(\tau_i,x)
\Biggr|^\beta \,\mathrm{d}x \quad
\mbox{and}\\
&&\int_{-\infty}^\infty\Biggl( \Biggl|\sum_{i=1}^k\theta_iL_\ast(\tau_i,x)
\Biggr|^\beta
\operatorname{sgn} \Biggl(\sum_{i=1}^k\theta_iL_\ast(\tau_i,x) \Biggr) \Biggr)\,\mathrm{d}x .
\end{eqnarray*}
\end{corollary}

\begin{pf}
If we let $ k_n:=n^{(1+\alpha)/\alpha} $, then for all $ n\in
\mathbb{N} $ and $ x\in\mathbb{Z} $, we have that
\[
\Gamma_n(\tau,x/n)=\int_0^\tau\mathbh{1}_{\{x/n\}}(X_n(t))\,\mathrm{d}t
=k_n^{-1}\int_0^{k_n\tau} \mathbh{1}_{\{x\}}(X(t))\,\mathrm{d}t
= n^{-(\alpha+1)/\alpha}\Gamma(k_n\tau,\{x\}) .
\]
The result then follows from Proposition \ref{PrinceProp}.
\end{pf}

\section{The finite-dimensional distributions}

In this section, we prove the convergence of the finite-dimensional
distributions of $ \Xi_n $ toward the
finite-dimensional distributions of $ \Xi_\ast$. In order to do so,
we first compute the exact expression of
the finite-dimensional distributions of $ \Xi_\ast$.
The proofs in this section follow the ideas given in Kesten and Spitzer (\citeyear{KesSpi1979}).

In the \hyperref[sec1]{Introduction}, we defined
\[
\Xi_\ast(\tau):=\int_0^\infty L_\ast(\tau,x-)\,\mathrm{d}Z_+(x)+\int
_0^\infty L_\ast(\tau,-(x-))\,\mathrm{d}Z_-(x) ,
\]
where $ \{Z_+(t);t\geq0\} $ and $ \{Z_-(t);t\geq0\} $ are independent
copies of the $ \beta$-stable L\'{e}vy process, which can be
associated with the stable distribution $ \vartheta_\beta$ with
characteristic function given by
\[
\psi(\theta)=\exp\bigl(-|\theta|^\beta\bigl(A_1+\mathrm{i}A_2\operatorname{sgn}(\theta)\bigr) \bigr) .
\]

\begin{lemma} \label{FinitDistriLem}
For $ t_1,\ldots,t_k\geq0 $ and $ \theta_1,\ldots,\theta_k\in
\mathbb{R} $, we have that
\begin{eqnarray*}
&&\mathbb{E} \Biggl[\exp\Biggl(\mathrm{i}\sum_{j=1}^k\theta_j\Xi_\ast(t_j) \Biggr)
\Biggr]\\
&&\quad= \mathbb{E} \Biggl[\exp\Biggl(-A_1\int_{-\infty}^\infty\Biggl|\sum_{j=1}^k\theta_jL_\ast
(t_j,x) \Biggr|^\beta \,\mathrm{d}x \Biggr)\\
&&\qquad\hphantom{\mathbb{E} \Biggl[}{}\times\exp\Biggl(-\mathrm{i}A_2\int_{-\infty}^\infty\Biggl|\sum_{j=1}^k\theta_jL_\ast
(t_j,x) \Biggr|^\beta \,\mathrm{d}x
\operatorname{sgn} \Biggl(\sum_{j=1}^k\theta_jL_\ast(t_j,x) \Biggr) \Biggr) \Biggr] .
\end{eqnarray*}
\end{lemma}

\begin{pf}
The proof is similar to that given in Kesten and Spitzer (\citeyear{KesSpi1979}) (see
page 16ff).
Let $ \nu$ be the L\'{e}vy measure of $ Z_+ $. One can truncate the L\'{e}vy measure as follows:
\[
\nu_1(B)=\nu(B\cap\{y\in\mathbb{R};|y|\leq1\}) \quad \mbox{and} \quad
\nu_2(B)=\nu(B\cap\{y\in\mathbb{R};|y|>1\}).
\]
Let $ M(t) $ and $ A(t) $ be independent L\'{e}vy processes, with
respective characteristic functions
\[
\mathbb{E} \bigl[\mathrm{e}^{\mathrm{i}\theta M(t)} \bigr]
=\exp\biggl(t\int_{|y|\leq1} (\mathrm{e}^{\mathrm{i}\theta y}-1-\mathrm{i}\theta y )\nu_1(\mathrm{d}y) \biggr)
\]
and
\[
\mathbb{E} \bigl[\mathrm{e}^{\mathrm{i}\theta A(t)} \bigr]=\exp\biggl(t\int_{|y|\leq1} (\mathrm{e}^{\mathrm{i}\theta
y}-1 )\nu_2(\mathrm{d}y) \biggr),
\]
such that
\[
Z^+(t)=M(t)+A(t)+Dt ,
\]
where $ D $ is a suitable real constant. This decomposition exists and
is called the L\'{e}vy--It\^{o}
representation of $ Z^+ $.
The advantage of this representation is that $ M(t) $ is a martingale
and has all moments and $ A(t) $ is
a process with bounded variation. Since the process $ \{L_\ast
(t,x-);x\geq0\} $ is left-continuous
and independent with respect to the filtration $ \mathcal{F}_t $ generated
by $ Z^+(t) $, the process
$ \{L_\ast(t,x-);x\geq0\} $ is $ \mathcal{F}_t $-predictable. Moreover,
$ \{L_\ast(t,x-);x\geq0\} $ has bounded
support $ \mathbb{P} $-almost surely. Therefore, we can find a suitable
sequence of partitions
$ \{x_l^{(n)};l\in\mathbb{N}\}$, $n\in\mathbb{N} $, with $
x^{(n)}_l<x^{(n)}_{l+1} $ for all
$ l,n\in\mathbb{N} $ satisfying
\[
\lim_{l\rightarrow\infty} x_l^{(n)}=\infty\quad \mbox{and}\quad
\lim_{n\rightarrow\infty}\max_{l\in\mathbb{N}} \bigl(x_{l+1}^{(n)}-x_l^{(n)} \bigr)=0
\]
such that
\[
\int_0^\infty L_\ast(t,x-)\,\mathrm{d}M(x)=\lim_{n\rightarrow\infty}
\sum_{l=1}^\infty L_\ast\bigl(t,x_l^{(n)}-\bigr) \bigl(M\bigl(x_{l+1}^{(n)}\bigr)-M\bigl(x_l^{(n)}\bigr)\bigr)
\]
with probability 1
(see Meyer (\citeyear{Mey1976}), Chapter II, Section 23).
Moreover, we can also assume that
\[
\int_0^\infty L_\ast(t,x-)\,\mathrm{d}A(x)=\lim_{n\rightarrow\infty}
\sum_{l=1}^\infty L_\ast\bigl(t,x_l^{(n)}-\bigr) \bigl(A\bigl(x_{l+1}^{(n)}\bigr)-A\bigl(x_l^{(n)}\bigr)
\bigr)
\]
with probability 1.

From those considerations, it follows that there exists a sequence of partitions
$ (x_l^{(n)})_{l\in\mathbb{N}} $ such that
\[
\int_0^\infty L_\ast(t,x-)\,\mathrm{d}Z_+(x)=\lim_{n\rightarrow\infty}
\sum_{l=1}^\infty L_\ast\bigl(t,x_l^{(n)}-\bigr)
\bigl(Z_+\bigl(x_{l+1}^{(n)}\bigr)-Z_+\bigl(x_l^{(n)}\bigr) \bigr)
\]
with probability 1.
Since the increments $ D^{(n)}_l:=Z_+(x_{l+1}^{(n)})-Z_+(x_l^{(n)}),\
l\in\mathbb{N}, $ are independent and
have characteristic function
\[
\mathbb{E} \bigl[\mathrm{e}^{\mathrm{i}\theta D^{(n)}_l} \bigr]
=\exp\bigl(-\bigl(x_{l+1}^{(n)}-x_l^{(n)}\bigr)|\theta|^\beta\bigl(A_1+\mathrm{i}A_2\cdot\operatorname{sgn}(\theta)\bigr) \bigr)
\]
by dominated convergence, we have
\begin{eqnarray*}
&& \mathbb{E} \Biggl[\exp\Biggl(\mathrm{i}\sum_{j=1}^k\theta_j\int_0^\infty L_\ast
(t_j,x-)\,\mathrm{d}Z_+(x) \Biggr) \Biggr]\\
&&\quad=\lim_{n\rightarrow\infty}
\mathbb{E} \Biggl[\exp\Biggl(
\sum_{l=1}^\infty\sum_{j=1}^k\mathrm{i}\theta_jL_\ast\bigl(t_j,x_l^{(n)}-\bigr)
\bigl(Z_+\bigl(x_{l+1}^{(n)}\bigr)-Z_+\bigl(x_l^{(n)}\bigr) \bigr) \Biggr) \Biggr]\\
&&\quad=\lim_{n\rightarrow\infty}
\mathbb{E} \Biggl[\exp\Biggl(-\sum_{l=1}^\infty\bigl(x_{l+1}^{(n)}-x_l^{(n)} \bigr)
\Biggl|\sum_{j=1}^k\theta_jL_\ast\bigl(t_j,x_l^{(n)}-\bigr) \Biggr|^\beta\\
&&\qquad\hphantom{\lim_{n\rightarrow\infty}
\mathbb{E} \Biggl[\exp\Biggl(-\sum_{l=1}^\infty}
{}\times\Biggl(A_1+\mathrm{i}A_2\cdot\operatorname{sgn} \Biggl(\sum_{j=1}^k\theta_jL_\ast\bigl(t_j,x_l^{(n)}-\bigr)
\Biggr) \Biggr) \Biggr) \Biggr].\\
&&\quad= \mathbb{E} \Biggl[\exp\Biggl(-A_1\int_0^\infty\Biggl|
\sum_{j=1}^k\theta_jL_\ast\bigl(t_j,x_l^{(n)}\bigr) \Biggr|^\beta \,\mathrm{d}x\\
&&\qquad\hphantom{\mathbb{E} \Biggl[\exp\Biggl(}
{}-\mathrm{i}A_2\int_0^\infty\Biggl|\sum_{j=1}^k\theta_jL_\ast\bigl(t_j,x_l^{(n)}\bigr)
\Biggl|^\beta
\operatorname{sgn} \Biggl(\sum_{j=1}^k\theta_jL_\ast\bigl(t_j,x_l^{(n)}\bigr) \Biggr)\,\mathrm{d}x \Biggr) \Biggr].
\end{eqnarray*}
For $ Z_- $, one can proceed with similar arguments.
\end{pf}

\begin{proposition} \label{FinitDistriConvProp}
The finite-dimensional distributions of the processes $ \{\Xi
_n(t);t\geq0\} $ converge toward the finite-dimensional distributions
of the process $ \{\Xi_\ast(t);t\geq0\} $.
\end{proposition}

\begin{pf}
As in the previous sections, we define $ k_n:=n^{(1+\alpha)/\alpha} $ and
$ \kappa:=\frac{1}{\alpha}+\frac{1}{\beta} $.
We already saw that we can use the occupation time $ \{\Gamma(t,\{x\}
);t\geq0,x\in\mathbb{R}\} $ of the
process $ \{X(t);t\geq0\} $ to represent the process $ \{\Xi(t);t\geq
0\} $ as follows:
\[
\Xi(t)=\sum_{x\in\mathbb{Z}}\Gamma(t,\{x\})\xi(x) .
\]
It follows that
\[
\Xi_n(t)=n^{-\kappa}\Xi(k_nt)=n^{-\kappa}\sum_{x\in\mathbb
{Z}}\Gamma(k_nt,\{x\})\xi(x) .
\]
Let $ \varphi(\theta):=\mathbb{E} [\exp(\mathrm{i}\theta\xi(1)) ] $ be the
characteristic function of the scenery
random variable $ \xi(1) $.
It then follows from the above representation that
\[
\sum_{j=1}^k\theta_j\Xi_n(t_j)
=n^{-\kappa}\sum_{x\in\mathbb{Z}}\sum_{j=1}^k\theta_j\Gamma
(k_nt_j,\{x\})\xi(x)
\]
and
\begin{eqnarray*}
R_n:= \mathbb{E} \Biggl[\exp\Biggl(\mathrm{i}\sum_{j=1}^k\theta_j\Xi_n(t_j) \Biggr) \Biggr]
=\mathbb{E} \Biggl[\prod_{x\in\mathbb{Z}}\varphi\Biggl(n^{-\kappa}\sum
_{j=1}^k\theta
_j\Gamma(k_nt_j,\{x\}) \Biggr) \Biggr].
\end{eqnarray*}
The random scenery $ \{\xi(z);z\in\mathbb{Z}\} $ is in the domain of
attraction of a $ \beta$-stable
distribution with characteristic function given by
\[
\psi(\theta)=\exp\bigl(-|\theta|^\beta\bigl(A_1+\mathrm{i}A_2\cdot\operatorname{sgn}(\theta)
\bigr)\bigr) .
\]
This implies that
\[
1-\varphi(\theta)\sim|\theta|^\beta
\bigl(A_1+\mathrm{i}A_2\cdot\operatorname{sgn}(\theta)\bigr)\qquad
\mbox{as } \theta\rightarrow0 .
\]
Thus
\begin{eqnarray*}
\log(\varphi(\theta)) \sim\log(\psi(\theta)) \qquad\mbox{as }
\theta\rightarrow0 .
\end{eqnarray*}
Therefore, for $ |\theta|\leq1 $, we have that
\begin{eqnarray*}
\biggl|\frac{\log(\varphi(\theta))-\log(\psi(\theta))}{\log(\psi
(\theta))} \biggr|=\mathrm{o}(\theta) .
\end{eqnarray*}
If we define
\[
\varphi_{x,n}:=\varphi\Biggl( n^{-\kappa}\sum_{j=1}^k\theta_j\Gamma
(k_nt_j,\{x\}) \Biggr)
\]
and
\[
\psi_{x,n}:=\exp\Biggl(-n^{-\kappa\beta} \Biggl|\sum_{j=1}^k\theta_j\Gamma
(k_nt_j,\{x\}) \Biggr|^\beta
\Biggl(A_1+\mathrm{i}A_2\cdot\operatorname{sgn} \Biggl(\sum_{j=1}^k\theta_j\Gamma(k_nt_j,\{x\}) \Biggr)
\Biggr) \Biggr)
\]
for all $ x\in\mathbb{Z} $, one has
\begin{eqnarray*}
\biggl|\frac{\log(\varphi_{x,n})-\log(\psi_{x,n})}{\log(\psi_{x,n})} \biggr|
=\mathrm{o} \Biggl(n^{-\kappa}\sum_{j=1}^k\theta_j\Gamma(k_nt_j,\{x\}) \Biggr) .
\end{eqnarray*}
This implies
that
\begin{eqnarray*}
\biggl|\log\biggl(\prod_{x\in\mathbb{Z}}\varphi_{x,n} \biggr)
-\log\biggl(\prod_{x\in\mathbb{Z}}\psi_{x,n} \biggr) \biggr|
&=& \biggl|\sum_{x\in\mathbb{Z}}\log(\varphi_{x,n})-\sum_{x\in\mathbb
{Z}}\log(\psi_{x,n}) \biggr| \\
&\leq&\sum_{x\in\mathbb{Z}}\log(\psi_{x,n})
\mathrm{o} \Biggl(n^{-\kappa}\sum_{j=1}^k\theta_j\Gamma(k_nt_j,\{x\}) \Biggr).
\end{eqnarray*}
By Corollary \ref{PrinceKor}, the right-hand side of the previous
inequality converges toward zero in probability.
The continuity of the logarithm then implies that
\begin{eqnarray*}
\biggl|\prod_{x\in\mathbb{Z}}\varphi_{x,n}-\prod_{x\in\mathbb{Z}}\psi
_{x,n} \biggr|\longrightarrow0
\qquad \mbox{in probability as } n\rightarrow\infty.
\end{eqnarray*}
We use this and dominated convergence to prove that the limit of the
sequence $ \{R_n;n\in\mathbb{N}\} $
exists and is equal to the limit of the sequence
\begin{eqnarray*}
Q_n:=\mathbb{E} \Biggl[\exp\Biggl(-\sum_{x\in\mathbb{Z}}n^{-\kappa\beta} \Biggl|
\sum_{j=1}^k\theta_j\Gamma(k_nt_j,\{x\}) \Biggr|^\beta
\Biggl(A_1+\mathrm{i}A_2\cdot\operatorname{sgn} \Biggl(\sum_{j=1}^k\theta_j\Gamma(k_nt_j,\{x\}) \Biggr)
\Biggr) \Biggr) \Biggr].
\end{eqnarray*}
By Corollary \ref{PrinceKor} and Lemma \ref{FinitDistriLem}, the
sequence $ \{Q_n;n\in\mathbb{N}\} $
converges toward
\begin{eqnarray*}
Q_\ast&:=&\mathbb{E} \Biggl[\exp\Biggl(-\int_{-\infty}^\infty\Biggl|\sum
_{j=1}^k\theta
_jL_\ast(t_j,x) \Biggr|^\beta
\Biggl(A_1+\mathrm{i}A_2\cdot\operatorname{sgn} \Biggl(\sum_{j=1}^k\theta_jL_\ast(t_j,x) \Biggr) \Biggr)\,\mathrm{d}x \Biggr)
\Biggr]\\
&=&\mathbb{E} \Biggl[\exp\Biggl(\mathrm{i}\sum_{j=1}^k\theta_j\Xi_\ast(t_j) \Biggr) \Biggr].
\end{eqnarray*}
As we have seen in Lemma \ref{FinitDistriLem}, $ Q_\ast$ is the
characteristic function for the
finite-dimensional distributions of $ \{\Xi_\ast(t);t\geq0\} $. This
completes the proof of the proposition.
\end{pf}

\section{The tightness}

In this section, we prove that the sequence $ \{\Xi_n(t);t\geq0\} $
is tight.
The proof of Theorem \ref{MT} then follows since we have already
obtained the convergence of the finite-dimensional distributions in the previous section.
The main proof of tightness also follows the ideas given in Kesten and Spitzer (\citeyear{KesSpi1979}).
We first need some suitable inequalities for the occupation times of $
X_\ast$.
However, the proofs of those inequalities differ from those given in
Kesten and Spitzer (\citeyear{KesSpi1979}).

\begin{lemma}\label{Lem2}
There exists a function $ \epsilon\dvtx\mathbb{R}^+\rightarrow\mathbb
{R}^+ $ with the properties
$ \epsilon(A)\rightarrow0 $ as $ A\rightarrow\infty$ and
\[
\mathbb{P} \bigl(\Gamma(s,\{x\})>0\ \mbox{for\ some}\ x\ \mbox{with}\ |x|>As^{
{\alpha
}/({1+\alpha})} \bigr)\leq\epsilon(A) \qquad \mbox{for all } s\geq0.
\]
\end{lemma}

\begin{pf}
For a positive real number $ x $, we denote by $ \lceil x\rceil$ the
smallest integer which is
greater or equal to $ x $. Obviously, for all $ s\geq0 $, we have
\begin{eqnarray*}
&& \mathbb{P} \bigl(\Gamma(s,\{x\})>0 \mbox{ for some } x\mbox{ with }
|x|>As^{\alpha/(1+\alpha)} \bigr) \\
&&\quad\leq \mathbb{P} \bigl(|X(r)|> As^{\alpha/(1+\alpha)}\mbox{ for
some }
r\leq s \bigr) \\
&&\quad\leq \mathbb{P} \bigl(|X(r)|> A \bigl( \bigl\lceil s^{\alpha/(1+\alpha)}
\bigr\rceil-1
\bigr)\mbox{ for some } r\leq\bigl\lceil s^{\alpha/(1+\alpha)} \bigr\rceil
^{(1+\alpha)/\alpha} \bigr) \\
&&\quad= \mathbb{P} \bigl( \bigl|X \bigl( \bigl\lceil s^{{\alpha}/({1+\alpha})} \bigr\rceil
^{(1+\alpha)/\alpha}u \bigr) \bigr|>
A \bigl\lceil s^{\alpha/(1+\alpha)} \bigr\rceil-A \mbox{ for some }
u\leq1 \bigr) \\
&&\quad\leq \mathbb{P} \Bigl(\sup_{r\leq1}\bigl|X_{n(s)}(r)\bigr|> A/2 \Bigr) \qquad
\mbox{for } s>1,
\end{eqnarray*}
with $ n(s):= \lceil s^{\alpha/(1+\alpha)} \rceil\rightarrow
\infty$ as $ s\rightarrow\infty$.
Since
\[
\mathbb{P} \Bigl(\sup_{r\leq1}|X_n(r)|> A/2 \Bigr)\longrightarrow\mathbb{P}
\Bigl(\sup_{r\leq
1}|X_\ast(r)|> A/2 \Bigr)\qquad \mbox{as }
n\rightarrow\infty,
\]
we can define
\[
\epsilon(A):=\sup_{s\geq0}\mathbb{P} \Bigl(\sup_{r\leq1}\bigl|X_{n(s)}(r)\bigr|>
A/2 \Bigr)\qquad
\mbox{for all } A>0 .
\]
This proves the statement of the lemma.
\end{pf}

\begin{lemma}\label{Lem3}
There exists a $ C>0 $ such that for all $ s\geq0 $, one has
\[
\sum_{x\in\mathbb{Z}}\mathbb{E} [\Gamma^2(s,\{x\}) ]\sim
Cs^{2-{\alpha}/({1+\alpha})} .
\]
\end{lemma}

\begin{pf}
For a positive real number $ x $, we denote by $ \lfloor x\rfloor$ its
integer part.
We know that for $ w(s):= \lfloor s^{\alpha/(\alpha+1)}
\rfloor$, one has
\begin{eqnarray*}
\frac{(w(s))^{2(\alpha+1)/\alpha}}{s^2}\sum_{x\in\mathbb
{Z}}\Gamma_{w(s)}^2\bigl(1,\{x/w(s)\}\bigr)
=s^{-2}\sum_{x\in\mathbb{Z}}\Gamma^2 \bigl((w(s))^{(\alpha+1)/{\alpha}},\{x\} \bigr)
\leq s^{-2}\sum_{x\in\mathbb{Z}}\Gamma^2(s,\{x\})
\end{eqnarray*}
and
\begin{eqnarray*}
s^{-2}\sum_{x\in\mathbb{Z}}\Gamma^2(s,\{x\})
& \leq& s^{-2}\sum_{x\in\mathbb{Z}}\Gamma^2 \bigl(\bigl(w(s)+1\bigr)^{
({\alpha+1})/{\alpha}},\{x\} \bigr)\\
&=& \frac{(w(s)+1)^{2(\alpha+1)/{\alpha}}}{s^2}\sum_{x\in
\mathbb{Z}}\Gamma_{w(s)+1}^2\bigl(1,\bigl\{x/\bigl(w(s)+1\bigr)\bigr\}\bigr) .
\end{eqnarray*}
Consequently, one has
\[
s^{-2}\sum_{x\in\mathbb{Z}}\mathbb{E} [\Gamma^2(s,\{x\}) ]
\sim\sum_{x\in\mathbb{Z}}\mathbb{E} \bigl[\Gamma_{w(s)}^2\bigl(1,\{x/w(s)\}
\bigr) \bigr]
= \sum_{x\in\mathbb{Z}}\mathbb{E} \bigl[\tilde{\Gamma}_{w(s)}^2\bigl(1,\{
x/w(s)\}\bigr) \bigr] .
\]
It follows from the layer cake representation and the remark after the
proof of Proposition
\ref{SignedCardTowardMeasureProp} that
\begin{eqnarray*}
w(s)\sum_{x\in\mathbb{Z}}\tilde{\Gamma}_{w(s)}^2\bigl(1,\{x/w(s)\}\bigr)
=\frac{1}{w(s)}\int_0^\infty\operatorname{card}
\bigl\{x\in\mathbb{Z}\dvtx
w^2(s)\tilde{\Gamma}^2_{w(s)}\bigl(1,\{x/w(s)\}\bigr)>c \bigr\}\,\mathrm{d}c
\end{eqnarray*}
converges $ \mathbb{P} \times\tilde{\mathbb{P} } $-almost surely toward
\[
\int_0^\infty\ell\bigl(x\in\mathbb{R}\dvtx\tilde{L}^2(1,x)>c \bigr)\,\mathrm{d}c
=\int_{\mathbb{R}}\tilde{L}_\ast^2(1,x)\,\mathrm{d}x.
\]
Dominated convergence and Fubini's theorem imply that
\[
w(s)\sum_{x\in\mathbb{Z}}\mathbb{E} \bigl[\tilde{\Gamma}_{w(s)}^2\bigl(1,\{
x/w(s)\}\bigr)
\bigr]\longrightarrow
\int_{\mathbb{R}}\mathbb{E} [\tilde{L}_\ast^2(1,x) ]\,\mathrm{d}x\qquad
\mbox{as }s\rightarrow\infty.
\]
Therefore,
\[
w(s)s^{-2}\sum_{x\in\mathbb{Z}}\mathbb{E} [\Gamma^2(s,\{x\})
]\longrightarrow
\int_{\mathbb{R}}\mathbb{E} [\tilde{L}_\ast^2(1,x) ]\,\mathrm{d}x\qquad \mbox{as }
s\rightarrow\infty.
\]
This proves the statement of the lemma.
\end{pf}

\begin{lemma} \label{Lem4}
\textup{(1)} For all $ \beta\in(0,2] $ and $ \rho>0 $, there exists a $ C_1>0
$ such that as $ n\rightarrow\infty$,
we have
\[
\bigl|\mathbb{E} \bigl[\xi(0)\mathbh{1}_{[-\rho,\rho]}
(n^{-1/\beta}\xi(0)) \bigr] \bigr|\sim
C_1n^{(1-\beta)/\beta} .
\]

\textup{(2)} For all $ \beta\in(0,2) $ and $ \rho>0 $, there exists a $ C_2>0
$ such that as $ n\rightarrow\infty$,
we have
\[
\bigl|\mathbb{E} \bigl[\xi^2(0)\mathbh{1}_{[-\rho,\rho]}
\bigl(n^{-{1/\beta}}\xi(0)\bigr) \bigr] \bigr|\sim
C_2n^{(2-\beta)/{\beta}}.
\]
\end{lemma}

\begin{pf}
The random variable $ \xi(0) $ is in the domain of attraction of a $
\beta$-stable random variable
with characteristic function given by
\[
\psi(\theta)=\exp\bigl(-|\theta|^\beta\bigl(A_1+\mathrm{i}A_2\operatorname{sgn}(\theta)\bigr)\bigr) ,
\]
with $ 0<A_1<\infty$ and $ |A_1^{-1}A_2|\leq\tan(\uppi\beta/2) $.
A consequence of this setting is that for $\beta>1 $, we have $
\mathbb{E}
[\xi(0)]=0 $.
Further, if $ \beta\in(0,2] $, then there exist $ B_1,B_2\geq0 $
such that
\[
\lim_{\rho\rightarrow\infty}\rho^\beta\mathbb{P} \bigl(\xi(0)\geq
\rho\bigr)= B_1\quad
\mbox{and}\quad
\lim_{\rho\rightarrow\infty}\rho^\beta\mathbb{P} \bigl(\xi(0)\leq
-\rho\bigr)= B_2 .
\]
For $ \beta=2 $, we have $ B_1=B_2=0$ since the decay of the tail
probabilities is exponential in that case.
For $ \beta\neq1 $, we then have that
\begin{eqnarray*}
\bigl|\mathbb{E} \bigl[\xi(0)\mathbh{1}_{[-\rho,\rho]}(n^{-
{1}/{\beta
}}\xi(0)) \bigr] \bigr|
&=& \int_0^{\rho n^{{1}/{\beta}}}\mathbb{P} \bigl(|\xi(0)|\geq
c\bigr)\,\mathrm{d}c\\
&\sim& (B_1+B_2)\int_0^{\rho n^{{1}/{\beta}}}c^{-\beta}\,\mathrm{d}c\\
&=& (B_1+B_2)(1-\beta)^{-1}\rho^{1-\beta}n^{({1}/{\beta
})(1-\beta)}.
\end{eqnarray*}
This proves the first statement for $ \beta\neq1 $. For $ \beta=1 $,
the statement is just our assumption
from the \hyperref[sec1]{Introduction}.

Moreover, by similar arguments for $ \beta\neq2 $, we have that
\begin{eqnarray*}
\bigl|\mathbb{E} \bigl[\xi^2(0)\mathbh{1}_{[-\rho,\rho
]}(n^{-
{1}/{\beta}}\xi(0)) \bigr] \bigr|
&\sim& (B_1+B_2)\int_0^{\rho n^{{1}/{\beta}}}c^{1-\beta}\,\mathrm{d}c\\
&=& (B_1+B_2)(2-\beta)^{-1}\rho^{2-\beta}n^{(1/\beta)
(2-\beta)}.
\end{eqnarray*}
This completes the proof of the second statement.
\end{pf}

\begin{proposition}
The distributions of the sequence $ \{\Xi_{n};n\in\mathbb{N}\} $ are
tight with respect to the
Skorohod topology.
\end{proposition}

\begin{pf}
We follow the method given in Kesten and Spitzer (\citeyear{KesSpi1979}). Let $ \epsilon
>0 $ be given.
By Lemma~\ref{Lem2}, there exists an $ A>0 $ such that
$ \epsilon(AT^{-\alpha/(1+\alpha)} )\leq\epsilon/4 $.
This implies that
\begin{eqnarray*}
&& \mathbb{P} \biggl(\Xi_n(t)\neq n^{-\kappa}\sum_{|x|\leq An}\Gamma
(k_nt,\{x\}
)\xi(x)\mbox{ for some } t\leq T \biggr) \\
&&\quad\leq \mathbb{P} \bigl(\Gamma(k_nT,\{x\})>0\mbox{ for some } x \mbox{ with }
|x|>Ak_n^{\alpha/(1+\alpha)} \bigr) \\
&&\quad\leq \epsilon\bigl(AT^{-\alpha/(1+\alpha)} \bigr)\\
&&\quad\leq \epsilon/4.
\end{eqnarray*}
There exists a $ \rho_0>0 $ with the property that for all $ \rho
>\rho_0 $ and all $ n\in\mathbb{N} $, we have
\[
3An\bigl(1-\mathbb{P} \bigl(-\rho n^{1/\beta}\leq\xi(0)\leq\rho
n^{1/\beta} \bigr)\bigr) \leq\epsilon/4 .
\]
This is valid since for suitable $ B_1,B_2\geq0 $, we have
\[
\lim_{\rho\rightarrow\infty}\rho^\beta\mathbb{P} \bigl(\xi(0)\geq
\rho\bigr)= B_1
\quad \mbox{and}\quad
\lim_{\rho\rightarrow\infty}\rho^\beta\mathbb{P} \bigl(\xi(0)\leq
-\rho\bigr)= B_2 .
\]
For all $ x\in\mathbb{Z} $, we have the random variables
\begin{eqnarray*}
\bar{\xi}_n(x)&:=&\xi(x)\mathbh{1}_{[-\rho,\rho
]}(n^{-1/\beta}\xi(x)) ,
\\
E_n&:=&n^{-\kappa}\frac{1}{T}\mathbb{E} \biggl[\sum_{x\in\mathbb
{Z}}\Gamma(k_nt,\{
x\})\bar{\xi}_n(x) \biggr]
=n^{-\kappa}\frac{1}{T}\mathbb{E} \biggl[\sum_{x\in\mathbb{Z}}
\Gamma(k_nt,\{x\})\mathbb{E} [\bar{\xi}_n(x) ] \biggr]
\end{eqnarray*}
and
\[
\bar{\Xi}_n(t):=n^{-\kappa}\sum_{x\in\mathbb{Z}}\Gamma(k_nt,\{x\})
\bigl(\bar{\xi}_n(x)-\mathbb{E} [\bar{\xi}_n(x) ] \bigr) .
\]

\textit{Claim} 1. The family of random variables $ \{E_n(t);n\in\mathbb
{N}\} $ is bounded. This is true
since, by Lemma \ref{Lem4}, we have
\begin{eqnarray*}
\biggl|\sum_{x\in\mathbb{Z}}\Gamma(k_nt,\{x\})\mathbb{E} [\bar{\xi }_n(x) ] \biggr|
&=& |\mathbb{E} [\bar{\xi}_n(0) ] |\sum_{x\in\mathbb{Z}}\Gamma
(k_nt,\{x\})\\
&=&k_nt |\mathbb{E} [\bar{\xi}_n(0) ] |
\leq Ctn^{(\alpha+1)/\alpha}n^{(1/\beta)(1-\beta)}
\end{eqnarray*}
and
$ \frac{\alpha+1}{\alpha}+\frac{1}{\beta}(1-\beta)-\kappa=0.$\vspace*{1pt}

\textit{Claim} 2. For all $ \eta>0 $, there exists an $ n_0\in\mathbb
{N} $ such that for all $ n\geq n_0 $, we have
\[
\mathbb{P} \biggl(\sup_{t\leq T}|\Xi_n(t)-\bar{\Xi}_n(t)-E_nt|>\frac
{\eta}{2}
\biggr)\leq\frac{\epsilon}{2} .
\]
To see this, we first note that
\[
\Xi_n(t)-\bar{\Xi}_n(t)-E_nt = n^{-\kappa}\sum_{x\in\mathbb
{Z}}\Gamma(k_nt,\{x\})
\bigl(\xi(x)-\bar{\xi}_n(x) \bigr)
\]
since
\begin{eqnarray*}
&& \Xi_n(t)-\bar{\Xi}_n(t)-E_nt-n^{-\kappa}\sum_{x\in\mathbb
{Z}}\Gamma(k_nt,\{x\})
\bigl(\xi(x)-\bar{\xi}_n(x) \bigr) \\
&&\quad= n^{-\kappa} \biggl(\sum_{x\in\mathbb{Z}}\Gamma(k_nt,\{x\})\mathbb
{E} [\bar
{\xi}(x) ]
-\frac{t}{T}\mathbb{E} \biggl[\sum_{x\in\mathbb{Z}}\Gamma(k_nt,\{x\}
)\mathbb{E} [\bar
{\xi}(x) ] \biggr] \biggr)\\
&&\quad= n^{-\kappa}\mathbb{E} [\bar{\xi}(0) ] \biggl(\sum_{x\in\mathbb
{Z}}\Gamma
(k_nt,\{x\})
-\frac{t}{T}\mathbb{E} \biggl[\sum_{x\in\mathbb{Z}}\Gamma(k_nt,\{x\}) \biggr]
\biggr)\\
&&\quad= n^{-\kappa}\mathbb{E} [\bar{\xi}(0) ] \biggl(k_nt-\frac{t}{T}k_nT \biggr)\\
&&\quad=0.
\end{eqnarray*}
Lemma \ref{Lem4} implies that
\begin{eqnarray*}
&& \mathbb{P} \biggl(n^{-\kappa}\sum_{x\in\mathbb{Z}}\Gamma(k_nt,\{x\})
\bigl(\xi
(x)-\bar{\xi}_n(x) \bigr)\neq0
\mbox{ for some } t\leq T \biggr) \\
&&\quad\leq \mathbb{P} \bigl(\Gamma(k_nT,\{x\})>0\mbox{ for some } x
\mbox{ with }|x|>Ak_n^{\alpha/(1+\alpha)} \bigr) \\
&&\qquad{}+ \mathbb{P} \bigl(\xi(x)\neq\bar{\xi}_n(x)\mbox{ for some }
|x|\leq Ak_n^{\alpha/(1+\alpha)} \bigr) \\
&&\quad\leq \epsilon\bigl(AT^{-\alpha/(1+\alpha)} \bigr)
+3Ak_n^{\alpha/(1+\alpha)}\mathbb{P} \bigl(\xi(0)\neq\bar{\xi
}_n(0) \bigr) \\
&&\quad\leq \frac{\epsilon}{4}+3An
\bigl(1-\mathbb{P} \bigl(-\rho n^{1/\beta}\leq\xi(0)\leq\rho n^
{1/\beta} \bigr) \bigr) \\
&&\quad\leq \frac{\epsilon}{2}.
\end{eqnarray*}

\textit{Claim} 3. There exists a $ K_0>0 $ such that for all $ n\in
\mathbb{N} $, we have
\[
\mathbb{E} [ |\bar{\Xi}_n(t_2)-\bar{\Xi}_n(t_1) |^2 ]\leq
C_0(t_2-t_1)^{2-({1+\alpha})/{\alpha}} .
\]
We define the $ \sigma$-field $ \mathcal{X}=\{X(t);t\geq0\} $. It then
follows from
the independence of $ \{X(t);t\geq0\} $ and $ \{\xi(x);x\in\mathbb
{Z}\} $ that
\begin{eqnarray*}
&& \mathbb{E} \biggl[ \biggl(\sum_{x\in\mathbb{Z}}\bigl(\Gamma(k_nt_2,\{x\})-\Gamma
(k_nt_1,\{x\})\bigr)\bar{\xi}_n(x) \biggr)^2 \biggr]\\
&&\quad=\mathbb{E} \biggl[\mathbb{E} \biggl[ \biggl(
\sum_{x\in\mathbb{Z}}\bigl(\Gamma(k_nt_2,\{x\})
-\Gamma(k_nt_1,\{x\})\bigr)\bar{\xi}_n(x) \biggr)^2 \bigg|\mathcal{X} \biggr] \biggr]\\
&&\quad= \mathbb{E} \biggl[\sum_{x\in\mathbb{Z}}\bigl(\Gamma(k_nt_2,\{x\})
-\Gamma(k_nt_1,\{x\})\bigr)^2\mathbb{E} [\bar{\xi}^2_n(x) |\mathcal{X} ] \biggr]\\
&&\quad= \sum_{x\in\mathbb{Z}} \mathbb{E} \bigl[\bigl(\Gamma(k_nt_2,\{x\})
-\Gamma(k_nt_1,\{x\})\bigr)^2 \bigr]\mathbb{E} [\bar{\xi}^2_n(x) ].\\
\end{eqnarray*}
This implies that
\begin{eqnarray*}
\mathbb{E} [ |\bar{\Xi}_n(t_2)-\bar{\Xi}_n(t_1) |^2 ]
&\leq& n^{-2\kappa}
\sum_{x\in\mathbb{Z}}\mathbb{E} \bigl[\bigl(\Gamma(k_nt_2,\{x\})-\Gamma
(k_nt_1,\{x\}
)\bigr)^2 \bigr]\mathbb{E} [\bar{\xi}_n^2(x) ]\\
&=& n^{-2\kappa}\mathbb{E} \biggl[\sum_{x\in\mathbb{Z}}\bigl(\Gamma(k_nt_2,\{
x\})
-\Gamma(k_nt_1,\{x\})\bigr)^2 \biggr]\mathbb{E} [\bar{\xi}_n^2(0) ] .
\end{eqnarray*}
Conditioned on $ \mathcal{A}:=\{\lambda_i;i\in\mathbb{Z}\} $, the
process $ X $ has the strong Markov property.
Using this, we can prove that for $ t_1\leq t_2 $, the conditional
distribution of
$ \sum_x(\Gamma(t_2,\{x\})-\Gamma(t_1,\{x\}))^2 $ with respect to $
\mathcal{A} $ equals the conditional distribution
of $ \sum_x \Gamma^2(t_2-t_1,\{x\}) $ with respect to $ \mathcal{A} $. Hence,
\begin{eqnarray*}
\mathbb{E} \biggl[\sum_{x\in\mathbb{Z}}\bigl(\Gamma(t_2,\{x\})-\Gamma(t_1,\{
x\})\bigr)^2 \biggr]
&=&\mathbb{E} \biggl[\mathbb{E} \biggl[\sum_{x\in\mathbb{Z}}\bigl(\Gamma(t_2,\{x\}
)-\Gamma(t_1,\{x\}
)\bigr)^2 \big|\mathcal{A} \biggr] \biggr]\\
&=& \mathbb{E} \biggl[\mathbb{E} \biggl[\sum_{x\in\mathbb{Z}}\Gamma
^2(t_2-t_1,\{x\}) \big|\mathcal{A}
\biggr] \biggr]\\
&=& \mathbb{E} \biggl[\sum_{x\in\mathbb{Z}}\Gamma^2(t_2-t_1,\{x\}) \biggr].
\end{eqnarray*}
By Lemma \ref{Lem3}, it follows that
\begin{eqnarray*}
\mathbb{E} \biggl[\sum_{x\in\mathbb{Z}}\bigl(\Gamma(k_nt_2,\{x\})-\Gamma
(k_nt_1,\{x\}
)\bigr)^2 \biggr]
&\leq& Ck_n^{2-\alpha/(1+\alpha)}(t_2-t_1)^{2-\alpha
/(1+\alpha)} \\
&=& Cn^{2(1+\alpha)/\alpha-1}(t_2-t_1)^{2-\alpha
/(1+\alpha)}.
\end{eqnarray*}
Moreover, we know that
\[
\mathbb{E} [\bar{\xi}_n^2(0) ]\leq\tilde{C}n^{(2-\beta)
(1/\beta)} .
\]
Putting this all together, we obtain
\[
\mathbb{E} [ |\bar{\Xi}_n(t_2)-\bar{\Xi}_n(t_1) |^2 ]
\leq C_0n^{(2-\beta)(1/\beta)}
n^{-2\kappa}n^{2(1+\alpha)/\alpha-1}(t_2-t_1)^{2-
{\alpha}/({1+\alpha})} .
\]
Since $ (2-\beta)\frac{1}{\beta}-2\kappa+2\frac{1+\alpha}{\alpha
}-1=0 $, Claim 3 follows.

Since $ 2-\frac{\alpha}{1+\alpha}> 1 $, the tightness in the
Skorohod topology of the family
$ \{\Xi_n;n\in\mathbb{N}\} $ now follows from Claims 1--3 and a
theorem of Billingsley (\citeyear{Bil1968})
(see page 95).
\end{pf}

\section*{Acknowledgements}
The authors wish to express their deepest gratitude toward the
probability group and the staff of
Academia Sinica and National Taiwan University for mathematical and
administrational help during their
visit to Taiwan. Special thanks go to Shieh Narn-Rueih, Hwang Chii-Ruey
and Sheu Shuenn-Jyi for many
interesting discussions on probability theory. Moreover, the authors
would like to thank the referee
for his very detailed report which helped to improve the manuscript.

\printhistory

\end{document}